\documentclass[a4paper,reqno, 11pt]{amsart}

\usepackage[margin=3cm]{geometry}
\usepackage{amssymb}
\usepackage{amsmath}
\usepackage{amsthm}
\usepackage{mathabx}
\usepackage{enumerate}
\usepackage{verbatim}
\usepackage{lmodern}

\usepackage{tikz}
\usepackage{tikz-cd}
\usetikzlibrary{matrix,calc,arrows,decorations.pathmorphing}
\tikzset{node distance=2cm, auto}

\tikzset{cd/.style=matrix of math nodes,row sep=2em,column sep=2em, text height=1.5ex, text depth=0.5ex}
\tikzset{cdar/.style=->,auto}
\tikzset{mid/.style={anchor=mid}} 
\tikzset{narrowfill/.style={inner sep=1pt, fill=white}}

\newtheorem{theorem}{Theorem}[section]
\newtheorem{lemma}[theorem]{Lemma}
\newtheorem{corollary}[theorem]{Corollary}

\newtheorem{proposition}[theorem]{Proposition}
\theoremstyle{definition}
\newtheorem{definition}[theorem]{Definition}
\newtheorem{example}[theorem]{Example}

\newcommand{\RR}{\mathbb{R}}
\newcommand{\CC}{\mathbb{C}}
\newcommand{\NN}{\mathbb{N}}
\newcommand{\ZZ}{\mathbb{Z}}

\newcommand{\Hh}{\mathcal{H}}
\newcommand{\Kk}{\mathcal{K}}

\newcommand{\Aa}{\mathcal{A}}

\newcommand{\Ll}{\mathcal{L}}

\newcommand{\Bb}{\mathcal{B}}

\newcommand{\Ss}{\mathcal{S}}
\newcommand{\Mm}{\mathcal{M}}
\newcommand{\Nn}{\mathcal{N}}
\newcommand{\Ee}{\mathcal{E}}
\newcommand{\Ii}{\mathcal{I}}

\newcommand{\red}{\mathrm{r}}

\newcommand{\sspan}{\operatorname{span}\,}
\newcommand{\cspan}{\overline{\operatorname{span}}\,}

\newcommand{\eps}{\varepsilon}

\newcommand{\Id}{\operatorname{id}}

\newcommand{\Mloc}{M_{\mathrm{loc}}}
\newcommand{\loc}{\mathrm{loc}}
\newcommand{\mult}{\mathrm{mult}}
\newcommand{\Bis}{\mathrm{Bis}}

\newcommand{\ess}{\mathrm{ess}}

\newcommand{\alg}{\mathrm{alg}}

\newcommand{\Mod}{\mathrm{Mod}}
\newcommand{\llangle}{\langle\langle}
\newcommand{\rrangle}{\rangle\rangle}

\DeclareMathOperator{\clspan}{\overline{span}}

\theoremstyle{remark}
\newtheorem{remark}[theorem]{Remark}
\newtheorem{notation}[theorem]{Notation}

\usepackage[pdftitle={Aperiodic dynamical inclusions of $C^*$-algebras},
pdfauthor={Jonathan Taylor},
pdfsubject={Mathematics}]{hyperref}

\title{Aperiodic dynamical inclusions of $C^*$-algebras}

\author{Jonathan Taylor}
\address{Jonathan Taylor, Mathematisches Institut, Georg-August-Universität Göttingen, Bunsenstraße 3-5, 37073 Göttingen, Germany}
\email{jonathan.taylor@mathematik.uni-goettingen.de}

\thanks{The author was supported by a stipend from the German Academic Exchange Service (DAAD) Funding Program 57450037}

\begin{document}
	\maketitle

	\section*{Abstract}
	We define the local multiplier module of a Hilbert module in analogy to the local multiplier algebra for $C^*$-algebras. 
	We use properties of the local multiplier module to lift non-closed actions on $C^*$-algebras by Hilbert bimodules to closed actions on local multiplier algebras, and descend known results on such closed actions down to their unclosed counterparts. 
	We define aperiodic dynamical inclusions and characterise them as crossed products by inverse semigroup actions.
	We describe the slice structure for such inclusions and show that all slices arise as linear combinations of slices already present in the inverse semigroup action.
	
	\section{Introduction}
	
	One of the most common ways to construct interesting $C^*$-algebras is to build crossed products from dynamical systems.
	Given a pair of $C^*$-algebras $A\subseteq B$ one can ask whether $B$ `acts' in a meaningful way on $A$, and if so can $B$ be explicitly described in terms of these dynamics.
	With no further conditions on the $C^*$-algebras $A\subseteq B$, the answer is probably not. 
	For this reason, many authors have explored the question by adding several helper-conditions to the inclusion $A\subseteq B$ and in various capacities have answered the question positively.
	
	Two of these helper conditions have been so far inescapable. 
	The first is that $A$ should contain an approximate identity for $B$ (\emph{non-degeneracy}), and the second is that $B$ should be generated as a $C^*$-algebra by normalisers of $A$.
	A \emph{normaliser} of the inclusion $A\subseteq B$ is an element $m\in B$ such that $m^*Am$ and $mAm^*$ are contained in $A$.
	This ensures that one can consider closed $A$-submodules of normalisers, called \emph{slices}, without losing information about $B$ not encoded by these. 
	There are typically two more ingredients in the definition of these pairs, although these are the conditions which are modified most often in papers on the topic.
	One usually requires that $A$ is maximal in some way in $B$, and that there is a conditional expectation of some kind $P:B\to A$.
	The flavour of maximality changes from author to author, many papers building stronger results from the same conditions or showing similar results hold under weaker ones.
	The expectation is sometimes removed in favour of some other condition at surface level, but often one can build kinds of conditional expectations using the replaced condition.
	
	One of the first bids to describe these pairs of $C^*$-algebras was by Kumjian \cite{K1} and Renault \cite{R1} to study commutative Cartan subalgebras. 
	Renault considered pairs $A\subseteq B$ where $A$ is a maximal commutative subalgebra (\emph{masa}) in $B$, and there is a faithful conditional expectation $P:B\to A$. 
	In combination with results from \cite{K1}, Renault was able to describe these pairs in terms of twisted groupoid $C^*$-algebras for \'etale, locally compact, Hausdorff, second countable, effective groupoids.
	This was later expanded upon by Exel \cite{E1} to allow for noncommutative Cartan subalgebras by replacing (masa) with a virtual maximality condition and to describe them as reduced section algebras of Fell bundles over inverse semigroups.
	This was later expanded upon by Kwa\'sniewski and Meyer in \cite{KM1}, who showed that the Cartan pairs of Exel could be described as crossed products by closed and purely outer actions, as well as showing that Exel's virtual maximality condition is equivalent to a number of other conditions on the inclusion.
	
	Another condition one may put on the pair $A\subseteq B$ is \emph{aperiodicity}, which is explored in more detail in \cite{KM3}, \cite{KM2}, and \cite{KM4}. 
	This condition guarantees that there is at most one pseudo-expectation $E:B\to I(A)$, that is, a conditional expectation taking values in Hamana's injective hull of $A$.
	One particular subset of these expectations are expectations that take values in the local multiplier algebra $\Mloc(A)$ of $A$, which embeds in the injective hull of $A$ (cf. \cite{F1}). 
	
	In this article we aim to mimic results of the previously mentioned authors and papers while adjusting some of the conditions.
	We shall consider inclusions $A\subseteq B$ that are aperiodic and have a faithful conditional expectation taking values in the local multiplier algebra of $A$.
	One of the challenges that arises is that the conditional expectation \mbox{$E:B\to\Mloc(A)$} maps part of $B$ isomorphically onto a subalgebra of $\Mloc(A)$.
	This is undesireable since the inclusion $A\subseteq\Mloc(A)$ is quite badly behaved, as $\Mloc(A)$ does not have interesting dynamics on $A$.
	Thus we shall insist on a technical condition limiting the multiplicative domain of $E$ as described by Choi \cite{C1}. 
	This ensures that the `intersection' between $B$ and $\Mloc(A)$ is as small as possible, namely is exactly $A$ itself.
	
	One of the main tools we use throughout this article is an analogous construction of the local multiplier algebra applied to Hilbert bimodules.
	Given a Hilbert $A$--$B$-bimodule $X$, one can construct a Hilbert $\Mloc(A)$--$\Mloc(B)$-bimodule, called the \emph{local multiplier module of $X$}, analogously to the construction of the local multiplier algebra.
	We prove some useful properties of the local multiplier module of a Hilbert bimodule, following Ara and Matthiew \cite{AM1}.
	We also show how some bimodule-specific properties interact with that local multiplier module construction.
	
	Our main application of the local multiplier module construction is to take an aperiodic action of an inverse semigroup by Hilbert bimodules on a $C^*$-algebra $A$, and gain a corresponding Fell bundle over an inverse semigroup with unit fibre $\Mloc(A)$.
	To then gain an inverse semigroup action on $\Mloc(A)$ is not immediate, but does follow from an application of the main theorem of \cite{BE2}.
	We can then show that this induced action on $\Mloc(A)$ gives rise to an Exel-Cartan inclusion, and so we then have the results of \cite{E1}, \cite{KM1} which we can apply to this induced inclusion.
	We then show that with another technical condition (that the canonical conditional expectation has minimal multiplicative domain), much of the structure of the induced Cartan inclusion descends to the original action of interest.
	
	\subsection*{Acknowledgements}
	The author would ilke to thank his doctoral advisor Ralf Meyer for all the assistance and expertise he provided.
	This article consists of results from the author's PhD thesis \cite{T1}.
	
	\section{Local multiplier modules and actions}
	Analogous to the construction of the local multiplier algebra of a $C^*$-algebra (cf.\ \cite[Definition~2.3.1]{AM1}), we define the local multiplier module for Hilbert $C^*$-modules. 
	We recall some key properties of the local multiplier algebra both from \cite{AM1} and others that follow from short arguments. 
	We also briefly study the relationship between a $C^*$-algebra and its local multiplier algebra, and show how some of these properties translate into the Hilbert module setting.
	Additionally, we discuss some properties specific to Hilbert (bi)modules and how these affect properties of local multiplier modules.
	
	\subsection{Local multiplier algebras}
	
	Throughout, $A$ will denote a fixed $C^*$-algebra, and $\Ii_e(A)$ will be the lattice of essential ideals of $A$ ordered by containment.
	Note that $\Ii_e(A)$ is directed, as the intersection of two essential ideals is itself an essential ideal.
	
	\begin{definition}[{\cite[Definition~2.3.1]{AM1}}]
		Let $I,J\in\Ii_e(A)$ with $J\subseteq I$.
		The restriction map $M(I)\to M(J)$ is injective and the collection of these restrictions over the net of essential ideals gives rise to the inductive limit
		$$\Mloc(A):=\varinjlim\limits_{I\in\Ii_e(A)}M(I),$$
		which we call the \emph{local multiplier algebra of $A$}.
	\end{definition}
	
	By definition, the multiplier algebras $M(I)$ canonically embed into $\Mloc(A)$ for all essential ideals $I\triangleleft A$. Moreover, if $J\triangleleft A$ is not essential, the ideal $J\oplus J^\perp$ is essential and we have $M(J)\subseteq M(J)\oplus M(J^\perp)=M(J\oplus J^\perp)$, which then includes into $\Mloc(A)$. In this fashion, one may think of $\Mloc(A)$ as the $C^*$-algebra generated by multipliers on ideals of $A$. Considering that ideals of $A$ are in bijective correspondence with open subsets of the spectrum $\hat{A}$ of $A$, we see that $\Mloc(A)$ is densely spanned by multipliers defined on these open subsets, hence the name \emph{local}. 

	\begin{lemma}[{\cite[Lemma~2.3.2]{AM1}}]\label{lem:AdetectsIdealsInMloc}
		Let $\tilde{I}\triangleleft\Mloc(A)$ be an ideal. Then $\tilde{I}\cap A=\{0\}$ if and only if $\tilde{I}=0$. Moreover, if $\tilde{I}$ is an essential ideal of $\Mloc(A)$, then $\tilde{I}\cap A$ is an essential ideal of $A$.
	\end{lemma}	
	
	\begin{lemma}[{\cite[Lemma~2.3.6]{AM1}}]\label{lem:locMultAddProd}
		For each $I\in\Ii_e(A)$ we have $\Mloc(I)=\Mloc(A)$. Let $(A_i)$ be a family of $C^*$-algebras. Then
		$$\Mloc\left(\bigoplus_{i}A_i\right)=\prod_{i}\Mloc(A_i).$$
	\end{lemma}
	
	Here the product $\prod_{i\in I} A_i$ of $C^*$-algebras $A_i$ consists of bounded families of elements $(a_i)_{i\in I}$, that is, each $a_i$ is an element of the $C^*$-algebra  $A_i$, and there is a constant $C\in\RR$ such that $||a_i||\leq C$ for all $i\in I$.
	This differs from the Cartesian product of the $A_i$ as sets, as that product may contain unbounded families.
	This boundedness criterium is needed to ensure that the product carries a norm, namely the supremum norm $||(a_i)_{i\in I}||:=\sup_{i\in I}||a_i||$.
	
	In the setting where $A=C_0(U)$ is a commutative $C^*$-algebra, essential ideals are of the form $C_0(V)$ for dense open subsets $V\subseteq U$, whereby the local multiplier algebra of $A$ loses the information of $\partial V$ for every dense open $V\subseteq U$.
	In the noncommutative setting, one may think of the quotient $A/I$ by an essential ideal $I\in\Ii_e(A)$ as the `boundary' of $I$.
	Lemma~\ref{lem:locMultAddProd} implies that even if $A/I$ is non-trivial, the quotient $\Mloc(A)/\Mloc(I)$ is always zero.
	
	\begin{lemma}\label{lem:AintersectMlocIdeal}
		Let $I\triangleleft A$ be an ideal. Then $A\cap\Mloc(I)=A\cap M(I)$, where the intersection is taken in $\Mloc(A)$.
		\begin{proof}
			Since $I\oplus I^\perp$ is an essential ideal in $A$ we have the inclusion $A\subseteq M(I\oplus I^\perp)=M(I)\oplus M(I)^\perp$ in $\Mloc(A)$. The left and right summands each embed respectively into the orthogonal summands $\Mloc(I)$ and $\Mloc(I^\perp)$, and so $A\cap\Mloc(I)\subseteq (M(I)\oplus M(I^\perp))\cap\Mloc(I)=M(I)$. Thus $A\cap\Mloc(I)=A\cap A\cap\Mloc(I)\subseteq A\cap M(I)$. The reverse inclusion holds as $M(I)\subseteq\Mloc(I)$.
		\end{proof}
	\end{lemma}	
	
	\begin{lemma}\label{lem:idealOfMlocContainment}
		Let $\tilde{I}\triangleleft\Mloc(A)$ be an ideal. Then $\tilde{I}\subseteq\Mloc(A\cap\tilde{I})$.
		\begin{proof}
			Lemma~\ref{lem:locMultAddProd} implies $\Mloc(A)=\Mloc(\tilde{I}\cap A)\oplus\Mloc((\tilde{I}\cap A)^\perp)$. Set $K:=\tilde{I}\cap\Mloc((\tilde{I}\cap A)^\perp)$. This is an ideal in $\Mloc(A)$ and so Lemma \ref{lem:AdetectsIdealsInMloc} gives $K\cap A=\{0\}$ if and only if $K=\{0\}$. Lemma~\ref{lem:AintersectMlocIdeal} then implies $K\cap A=\tilde{I}\cap\Mloc((\tilde{I}\cap A)^\perp)\cap A=(\tilde{I}\cap A) \cap M((\tilde{I}\cap A)^\perp)=\{0\}$. Thus $\tilde{I}$ must be contained in $\Mloc((\tilde{I}\cap A))$.
		\end{proof}
	\end{lemma}
	
	\subsection{Local multiplier modules}
	
	Let $A$ and $B$ be $C^*$-algebras. 
	If $X$ is a Hilbert $A$--$B$-bimodule we denote the right inner product by single angle brackets $\langle\cdot,\cdot\rangle$ and the left inner product by double angle brackets $\llangle\cdot,\cdot\rrangle$.
	We denote the source and range ideals $s(X)=\clspan\{\langle x,y\rangle:x,y\in X\}$ and $r(X)=\clspan\{\llangle x,y\rrangle:x,y\in X\}$.
	By definition of these ideals, the bimodule $X$ always induces a Morita equivalence between $s(X)$ and $r(X)$.
	
	The bidmodule structure also allows us to classify all subbimodules in terms of their range and source ideals.
	If $Y\subseteq X$ is a Hilbert $A$--$B$-subbimodule, one sees readily that $Y\subseteq r(Y)\cdot X$ and $Y\subseteq X\cdot s(Y)$.
	This is in fact equality, since for $x\in X$ and $y,z\in Y$ we have $x\cdot\langle y,z\rangle=\llangle x,y\rrangle\cdot z\in A\cdot Y\subseteq Y$, so $X\cdot s(Y)\subseteq Y$ (and similarly $r(Y)\cdot X\subseteq X$).
	
	For a Hilbert $A$--$B$-bimodule $X$, we denote by $X^*=\{x^*:x\in X\}$ its \emph{opposite module}.
	The opposite module $X^*$ has a $B$--$A$-bimodule structure.
	The addition in $X^*$ is given by $x^*+y^*=(x+y)^*$ for $x,y\in X$, the left $B$-action is given by $b\cdot x^*=(x\cdot b^*)^*$ for $b\in B$ and $x\in X$, and the right $A$-action is given by $x^*\cdot a=(a^*\cdot x)^*$ for $a\in A$ and $x\in X$.
	The left and right inner products are given by $\llangle x^*,y^*\rrangle=\langle x,y\rangle$ and $\langle x^*,y^*\rangle=\llangle x,y\rrangle$.
	The opposite bimodule satisfies $X\otimes_B X^*\cong r(X)$ and $X^*\otimes_A X\cong s(X)$, via maps defined on elementary tensors by $x\otimes y^*\mapsto \llangle x,y\rrangle$ and $x^*\otimes y\mapsto\langle x,y\rangle$ respectively for $x,y\in X$.
	
	One may always consider $A$ equipped with the canonical bimodule structure given by left and right multiplication, and inner products $\langle a,b\rangle=a^*b$ and $\llangle a,b\rrangle=ab^*$.
	The multiplier algebra $M(A)$ of $A$ can then be identified with the operators $A\to A$ that are adjointable with respect to the right inner product, when acting from the left (see \cite[Theorem~2.4]{Lance1}).
	For a more in depth introduction to Hilbert modules we recommend \cite{Lance1}.
	
	\begin{notation}
	Let $X$ and $Y$ be right Hilbert $B$-modules. 
	The rank-one operators $X\to Y$ are of the form $z\mapsto y\cdot\langle x,z\rangle$ for $x,z\in X$ and $y\in Y$.
	The space $\Kk(X,Y)$ of compact operators $X\to Y$ is the completion of the span of rank-one operators $X\to Y$, in the operator norm $||T||=\sup_{||x||=1}||Tx||$.
	We denote the space of adjointable operators $X\to Y$ by $\Ll(X,Y)$.
	If $X$ and $Y$ are Hilbert $A$--$B$-bimodules, and we wish to distinguish between operators that are left-adjointable and right-adjointable, we shall write
	\begin{align*}
		\Kk_R(X,Y)&=\{f:X\to Y| f \text{ is compact in the right Hilbert module structure}\}\\
		\Ll_R(X,Y)&=\{f:X\to Y| f \text{ is adjointable in the right Hilbert module structure}\}\\
		\Kk_L(X,Y)&=\{f:X\to Y| f \text{ is compact in the left Hilbert module structure}\}\\
		\Ll_L(X,Y)&=\{f:X\to Y| f \text{ is adjointable in the left Hilbert module structure}\}.
	\end{align*}
	In the case where $X$ and $Y$ are bimodules, and no subscript of $L$ or $R$ is specified as in $\Kk(X,Y)$ and $\Ll(X,Y)$, we assume these refer to the spaces of maps compatible with the right Hilbert module structure of $X$ and $Y$.
	If $X=Y$, the spaces $\Kk_*(X)$ and $\Ll_*(X)$ for $*=L,R$ are $C^*$-algebras with multiplication given by composition, ${}^*$-operation given by taking adjoints, and operator norm.
	In this scenario we have $\Ll_L(X)\cong M(\Kk_L(X))\cong M(s(X))$ and $\Ll_R(X)\cong M(\Kk_R(X))\cong M(r(X))$ by \cite[Theorem~2.4]{Lance1}.
	\end{notation}
	
	Let $X$ be a Hilbert $A$--$B$-bimodule.
	The space $\Ll_R(s(X),X)$ carries a $M(r(X))$--$M(s(X))$-bimodule structure by defining
	$$(\tau\cdot T\cdot \sigma)a=\tau(T(\sigma a)),$$
	for $\tau\in M(r(X))$, $\sigma\in M(s(X))$, and $a\in s(X)$.
	For $S,T\in\Ll_R(s(X),X)$ we see that $S^*T$ is an adjointable operator on $s(X)$, so uniquely defines a multiplier $\langle S,T\rangle$ of $s(X)$.
	Similarly, $ST^*$ defines an adjointable operator on $X$ which in turn uniquely defines a multiplier on $r(X)\cong X\otimes_B X^*$, which we label $\llangle S,T\rrangle$.
	The assignments $(S,T)\mapsto \langle S,T\rangle$ and $(S,T)\mapsto\llangle S,T\rrangle$ define a Hilbert $M(r(X))$--$M(s(X))$-bimodule structure on $\Ll_R(s(X),X)$.
	
	We define a Hilbert $M(r(X))$--$M(s(X))$-bimodule structure on $\Ll_L(r(X),X)$ similarly, and this is eloquently denoted by expressing adjointable maps in $\Ll_L(r(X),X)$ in post-fix notation: for $P\in\Ll_L(r(X),X)$ and $a\in r(X)$, the map $P:r(X)\to X$ evaluated at $a$ is written $aP$.
	
	For $\tau\in M(r(X))$, $\sigma\in M(s(X))$, and $P\in\Ll_L(r(X),X)$ we define $\tau\cdot P\cdot \sigma$ as the map $a\mapsto((a\tau)P)\sigma$.
	This gives $\Ll_L(r(X),X)$ a $M(r(X))$--$M(s(X))$-bimodule structure.
	For $P,Q\in\Ll_L(r(X),X)$ the compositions $P^*Q$ and $PQ^*$ define unique multipliers $\langle P,Q\rangle\in M(s(X))$ and $\llangle P,Q\rrangle\in r(X)$ respectively, which give rise to a Hilbert $M(r(X))$--$M(s(X))$-bimodule structure on $\Ll_L(r(X),X)$.
	
	There are two isomorphisms $X\cong\Kk_R(s(X),X)$ and $X\cong\Kk_L(r(X),X)$ given by mapping $x\in X$ to the respective left and right creation operators of $x$.
	That is, $x\mapsto |x\rangle\in\Kk_R(s(X),X)$, where $|x\rangle$ is the function mapping $b\in s(X)$ to $|x\rangle b=x\cdot b$.
	Similarly $\llangle x|\in\Kk_L(r(X),X)$ is the operator mapping $a\in r(X)$ to $a\llangle x|:=a\cdot x$.
	Denote the inverses of these isomorphisms by $\Xi_R:\Kk_R(s(X),X)\to X$ and $\Xi_L:\Kk_L(r(X),X)\to X$ respectively.
	
	\begin{lemma}\label{lem-symmetricMultiplierModuleStructure}
		Let $X$ be a Hilbert $A$--$B$-bimodule.
		Then $\Ll_R(s(X),X)$ and $\Ll_L(r(X),X)$ are isometrically isomorphic as Hilbert $M(s(X))$--$M(r(X))$-bimodules via the map taking $T\in\Ll_R(s(X),X)$ to the operator $\check{T}:r(X)\to X$ given by
		$$r(X)\ni a\mapsto \Xi_R(aT)\in X.$$
		The inverse of this map is then given by mapping $P\in\Ll_L(r(X),X)$ to the map $\hat{P}\in\Ll_R(s(X),X)$ given by
		$$s(X)\ni b\mapsto\Xi_L(Pb)\in X.$$
		\begin{proof}
			Fix $T\in\Ll_R(s(X),X)$.
			Since $X$ and $\Kk_R(s(X),X)$ are isomorphic, they have equal range ideals.
			Thus the map $a\mapsto\Xi_R(aT)$ is well defined since $aT\in r(X)\cdot\Ll_R(s(X),X)=\Kk_R(s(X),X)$.
			To see that $\check{T}$ is adjointable, fix $x\in X$ and $a\in r(X)$.
			Then
			$$\llangle a\check{T},x\rrangle=\llangle\Xi_R(aT),x\rrangle=\llangle aT,|x\rangle\rrangle,$$
			which is the unique multiplier on $r(X)$ specified by $aT |x\rangle^*=a(|x\rangle\circ T^*)^*$ (note that this is indeed a multiplier, since we can express it as a composition of adjointable maps).
			Thus $\check{T}$ has adjoint $X\to r(X)$ given by mapping $x$ to the unique multiplier of $r(X)$ associated to $|x\rangle\circ T^*$.
			This belongs to $r(X)=r(X)^*r(X)=\Kk(r(X))$ and not solely $M(r(X))$ since $|x\rangle$ is a compact operator, so the composition is also.
			
			For $S,T\in\Ll_R(s(X),X)$ and $x\in X$ we have $x\langle S^*,T\rangle=x(S^*\circ T)$ by definition, which is exactly $$\Xi_R((|x\rangle \circ S^*)\circ T)=\Xi_R(|x\rangle S^*T)=\Xi_R(|x\rangle)\langle S,T\rangle=x\langle S,T\rangle.$$
			Thus $S^*\circ T=\check{S}^*\circ\check{T}$ whereby $\langle S,T\rangle=\langle\check{S},\check{T}\rangle$ and the assignment $T\mapsto \check{T}$ is isometric.
			
			Lastly we show that $\check{\hat{P}}=P$ for all $P\in\Ll_L(r(X),X)$.
			This shall show that the map $T\mapsto\check{T}$ is a surjective map, hence an isomorphism, and the map $P\mapsto\hat{P}$ is its inverse.
			Fix $P\in\Ll_L(r(X),X)$.
			For all $a\in r(X)$ we have $a\check{\hat{P}}=\Xi_R(a\hat{P})$.
			This is the unique element of $X$ such that $\llangle \Xi_R(a\hat{P})|=a\check{\hat{P}}$.
			Hence $\Xi_R(a\hat{P})b$ is the unique element of $X$ satisfying $\Xi_R(a\hat{P})b=a\hat{P}b=a\Xi_L(Pb)=\Xi_L(aPb)$ for all $b\in s(X)$.
			Thus $a\check{\hat{P}}=\Xi_R(a\hat{P})=aP$ whereby $\check{\hat{P}}=P$.
		\end{proof}
	\end{lemma}
	
	Lemma~\ref{lem-symmetricMultiplierModuleStructure} allows us to consider $\Ll_R(s(X),X)$ and $\Ll_L(r(X),X)$ as an analogue of the multiplier algebra for Hilbert modules.
	
	If $K\subseteq J\triangleleft B$ are both essential ideals, then any adjointable operator $T\in\Ll(J,X\cdot J)$ restricts to a map $T|_K:K\to X\cdot J$.
	Moreover, $T|_K$ has range contained in $X\cdot K$ since each $a\in K$ can be written as $a=a_1a_2$ for some $a_1,a_2\in K$ and we have $T|_K a=(Ta_1)a_2\in X\cdot J\cdot K=X\cdot K$.
	The restriction $T|_K$ is adjointable from $K\to X\cdot K$, with adjoint $T^*|_{X\cdot K}$, which takes values in $K$ since $T^*(X\cdot K)=T^*(X\cdot J)\cdot K\subseteq K$.
	Thus, restriction of operators gives rise to a module homomorphism $\Ll(J,X\cdot J)\to\Ll(K,X\cdot K)$.
	If $T\in\Ll(J,X\cdot J)$ restricts to the zero operator on $K$, then $\{0\}=T^*(X\cdot K)=T^*(X\cdot J)\cdot K$, whereby $T^*(X\cdot J)$ is an ideal that annihilates $K$.
	But $K$ has zero annihilator in $J$ as $K$ is essential in $J$, whereby $T^*=0$.
	In particular, $T=0$ if and only if $T|_K=0$, and the map $\Ll(J,X\cdot J)\to\Ll(K,X\cdot K)$ is injective.
	These restriction maps give rise to an inductive system over essential ideals in $B$, ordered by reverse inclusion: $I\leq J$ if and only if $J\subseteq I$.
	
	\begin{definition}
		Let $X$ be a right Hilbert $B$-module. We define the \emph{(right-)local multiplier module} of $X$ as the inductive limit 
		$$X_\loc:=\varinjlim\limits_{J\in\Ii_e(B)}\Ll_R(J,X\cdot J)$$
		taken over essential ideals $J\triangleleft B$.
		
		Analogously, if $X$ is a left Hilbert $A$-module we define the \emph{(left-)local multiplier module} of $X$ as
		$${}_\loc X:=\varinjlim\limits_{I\in\Ii_e(A)}\Ll_L(I,I\cdot X).$$
	\end{definition}
	
	Taking the case $X=B$, with the canonical $B$-bimodule structure inherited from multiplication, we have that $\Ll_R(I,B\cdot I)=\Ll_R(I,I)=M(I)=\Ll_L(I,I)=\Ll_L(I,I\cdot B)$, thus $B_\loc={}_\loc B=\Mloc(B)$, and these constructions generalise the local multiplier algebra.
	
	\begin{lemma}\label{lem:starLocCommute}
		Let $X$ be a right Hilbert $B$-module. For each ideal $J\triangleleft B$ there is an isomorphism $\Ll_R(J,X\cdot J)^*\cong\Ll_L(J,J\cdot X^*)$ of left Hilbert $B$-modules. In particular $(X_\loc)^*\cong {}_\loc(X^*)$.
		\begin{proof}
			For $T\in\Ll_R(J,X\cdot J)$, define $S_T:J\to J\cdot X^*=(X\cdot J)^*$ by $S_Ta:=(Ta^*)^*$. Then the map $T\mapsto S_T$ gives an anti-isomorphism between $\Ll_R(J,X\cdot J)$ and $\Ll_L(J,J\cdot X^*)$.
			Moreover, these anti-isomorphisms entwine the restriction maps in the inductive limit, and so the inductive limits are isomorphic.
		\end{proof}
	\end{lemma}
	
	Many of the following results may be generalised to left Hilbert modules using either Lemma~\ref{lem:starLocCommute} or symmetry arguments.
	
	\begin{lemma}\label{lem:localisationCofinal}
		Let $X$ be a right Hilbert $B$-module. For any $J\in\Ii_e(B)$ we have $(X\cdot J)_\loc=X_\loc$.
		\begin{proof}
			This follows because the collection of essential ideals of $J$ is cofinal in $\Ii_e(A)$, over which the inductive limit is taken.
		\end{proof}
	\end{lemma}
	
	In the case where $X$ is a Hilbert bimodule, all subbimodules are of the form $X\cdot J$ for some ideal $J$. In particular, a subbimodule $Y\subseteq X$ has zero orthogonal complement in $X$ if and only if the source ideal of $Y$ is an essential ideal of the source of $X$. In the bimodule case, this then gives that the local multiplier module of any subbimodule with zero orthogonal complement is equal to the local multiplier module of the whole original module.
	
	\begin{corollary}\label{cor:locIdealDecomp}
		For any ideal $J\triangleleft B$ we have $X_\loc=(X\cdot (J\oplus J^\perp))_\loc$.
	\end{corollary}
	
	\begin{lemma}\label{lem:localisationAdditive}
		Let $X$ and $Y$ be right Hilbert $B$-modules. Then $(X\oplus Y)_\loc=X_\loc\oplus Y_\loc$.
		\begin{proof}
			This follows because the map $\Ll(J,X\cdot J)\oplus\Ll(J,Y\cdot J)\to \Ll(J,(X\oplus Y)\cdot J)$, $(T,S)\mapsto T+S$ is an isomorphism entwining the restriction maps for all ideals $J\triangleleft B$.
		\end{proof}
	\end{lemma}

	\begin{lemma}
		The module $X_\loc$ is a right Hilbert $\Mloc(B)$-module.
		\begin{proof}
			The argument at the start of this section shows that $\Ll_R(J,X\cdot J)$ is a right Hilbert $M(J)$-module for each essential ideal $J\triangleleft B$. 
			The inductive limit structure is preserved since all the inductive limit maps are restrictions of operators, which clearly preserve the right actions and inner products. Hence there is an induced right Hilbert $\Mloc(B)$-module structure on $X_\loc$.
		\end{proof}
	\end{lemma}
	
	In a Hilbert bimodule there are, at first glance, two ways to take an orthgonal complement of a subbimodule (one for each inner product).
	If $Y\subseteq X$ is a Hilbert subbimodule of a Hilbert bimodule $X$, and if $x\in X$ satisfies $\langle x,y\rangle=0$ for all $y\in Y$, then we also see that $\llangle y,x\rrangle\cdot z=y\langle x,z\rangle=0$ for all $y,z\in Y$.
	Thus $\llangle y,x\rrangle$ belongs to the annihilator of the range ideal $r(Y)$ for all $y\in Y$.
	However, the inner product $\llangle y,x\rrangle$ must also lie in $r(Y)$, as one may always write $y=ay'$ for some $a\in r(Y)$ and $y'\in Y$ by the Cohen-Hewitt factorisation theorem (see, for example, \cite[Lemma~4.4]{Lance1}).
	Thus if $x\in X$ annihilates $Y$ in the right inner product, it also does so in the left.
	Symmetrically, we see that orthogonal complements with respect to both inner products agree. 
	We shall denote both by $Y^\perp$.
	
	\begin{lemma}\label{lem:localPerpDecomp}
		Let $X$ be a Hilbert $A$--$B$-bimodule and let $Y\subseteq X$ be a closed Hilbert subbimodule. Then $(Y^\perp)_\loc=(Y_\loc)^\perp$ and $X_\loc=Y_\loc\oplus Y_\loc^\perp$.
		\begin{proof}
			We have $Y=X\cdot s(Y)$ and so by Lemma~\ref{lem:localisationAdditive} we gain 
			$$X_\loc=(X\cdot (s(Y)\oplus s(Y)^\perp))_\loc=Y_\loc\oplus (X\cdot s(Y)^\perp)_\loc.$$ 
			The right summand must then be equal to both $(Y^\perp)_\loc$ and $(Y_\loc)^\perp$.
		\end{proof}
	\end{lemma}
	
	\begin{lemma}\label{lem:locIsMultiplicativeForIdeals}
		Let $X$ be a Hilbert $A$--$B$-bimodule and let $J\triangleleft B$ be an ideal. Then \mbox{$(X\cdot J)_\loc= X_\loc\cdot\Mloc(J)$}.
		\begin{proof}
			First we show that $(X\cdot J)_\loc\subseteq X_\loc\cdot \Mloc(J)$. 
			Fix an essential ideal $J'\triangleleft B$, and define $J'_\ess:= JJ'\oplus (JJ')^\perp$, and note $J'_\ess$ is an essential ideal of $A$. 
			Fix $\xi\in\Ll(J'_\ess,(X\cdot J)\cdot J'_\ess)$. 
			Then $\xi=\xi|_{JJ'}\oplus\xi|_{(JJ')^\perp}$ and since $JJ'_\ess=JJ'$ we have $\xi|_{(JJ')^\perp}=0$, giving $\xi=\xi\cdot 1_J$. 
			Thus 
			$$\xi\in \Ll(J'_\ess,(X\cdot J)J'_\ess)\cdot\Mloc(J)\subseteq (X\cdot J)_\loc\cdot\Mloc(J).$$ 
			This gives $(X\cdot J)_\loc\subseteq X_\loc\cdot\Mloc(J)$.
			
			Lemmas~\ref{lem:localisationAdditive} and \ref{lem:localPerpDecomp} together imply that $X_\loc\cdot\Mloc(J)\oplus X_\loc\cdot\Mloc(J^\perp)=X_\loc=(X\cdot J)_\loc\oplus (X\cdot J^\perp)_\loc$. 
			For $\xi\in X_\loc$, $\tau\in\Mloc(J)$, and $\eta\in (X\cdot J^\perp)_\loc$ we have
			$$\langle\eta,\xi\cdot\tau\rangle=\langle\eta,\xi\rangle\tau\in\langle (X\cdot J^\perp)_\loc,X_\loc\rangle\cdot\Mloc(J).$$
			The set $\langle (X\cdot J^\perp)_\loc,X_\loc\rangle$ is contained in $\Mloc(J^\perp)$ by the above argument, and so the inner product $\langle\eta,\xi\cdot\tau\rangle$ is contained in $\Mloc(J^\perp)\cdot\Mloc(J)$, which is zero by Lemma~\ref{lem:locMultAddProd}.
			Thus $X_\loc\cdot\Mloc(J)\subseteq (X\cdot J)_\loc$ must hold.
		\end{proof}
	\end{lemma}

	\begin{lemma}
		Let $X$ be a Hilbert $A$--$B$-bimodule. Then $X_\loc$ and ${}_\loc X$ are Hilbert $\Mloc(A)$--$\Mloc(B)$-bimodules.
		\begin{proof}
			We show that $X_\loc$ is a left Hilbert $\Mloc(A)$-module. 
			The case for ${}_\loc X$ then follows by a symmetric argument, or from Lemma~\ref{lem:starLocCommute}. 
			For each essential ideal $J\triangleleft B$, let $I_J:=r(X\cdot J)\oplus r(X\cdot J)^\perp$. 
			Then $I_J$ is an essential ideal in $A$ and we claim $I_J\cdot X=X\cdot J$.
			To see this, note that $I_J$ satisfies $I_J\cdot X=r(X\cdot J)\cdot X\oplus r(X\cdot J)^\perp\cdot X$.
			The range ideal of $r(X\cdot J)^\perp\cdot X\cdot J$ is $r(X\cdot J)^\perp\cdot r(X\cdot J)=\{0\}$, so $r(X\cdot J)^\perp\cdot X\subseteq (X\cdot J)^\perp$.
			Since $J$ is essential in $s(X)$, we see $s((X\cdot J)^\perp)\cap J=\{0\}$ implies that $s((X\cdot J)^\perp)=\{0\}$, whereby $(X\cdot J)^\perp=\{0\}$ and $I_J\cdot X=X\cdot J$.
			
			For $\tau\in M(I_J)$ and $T\in\Ll_R(J,X\cdot J)$ we define $(\tau T)a:=\tau(Ta)$.
			This bilinear map commutes with the maps in the inductive system, giving the left action of $\Mloc(A)$ on $X_\loc$. 
			For $T,S\in\Ll_R(J,X\cdot J)$, we define the left inner product $\llangle T,S\rrangle$ as the element of $M(r(X\cdot J))$ corresponding to $TS^*$ under the isomorphism $\Ll(X\cdot J)\cong M(r(X\cdot J))$ arising from $X\cdot J\cong\Kk_R(J,X\cdot J)$.
			
			We also see clearly by construction that $\llangle T,S\rrangle P=TS^*P=T\langle S,P\rangle$ for all $T,S,P\in X_\loc$, giving the necessary compatibility of left and right inner products.
		\end{proof}
	\end{lemma}

	We now show that if $X$ is a Hilbert bimodule, then the left and right local multiplier modules agree.
	
	\begin{proposition}\label{prop:leftRightLocAgree}
		Let $X$ be a Hilbert $A$--$B$-bimodule. 
		For each ideal $J\triangleleft s(X)$ and the corresponding ideal $I_J:=r(X\cdot J)$, the isomorphisms $\Ll_R(J,X\cdot J)\cong \Ll_L(I_J,I_J\cdot X)$ of Hilbert $M(I_J)$--$M(J)$-bimodules from Lemma~\textup{\ref{lem-symmetricMultiplierModuleStructure}} preserve the inductive limit structures of $X_\loc$ and ${}_\loc X$. 
		In particular, ${}_\loc X\cong X_\loc$.
		\begin{proof}
			Let $\Ll_R(J,X\cdot J)\ni T\mapsto\check{T}\in\Ll_L(I_J,I_J\cdot X)$ denote the isomorphism from Lemma~\ref{lem-symmetricMultiplierModuleStructure}.
			
			Fix an essential ideal $J'\triangleleft J$ and note $I'=r(X\cdot J')$ is essential in $I$.
			If $T\in\Ll_R(J,X\cdot J)$ then $T|_{J'}$ gives the operator $\widecheck{T|_{J'}}\in\Ll_L(I',I'\cdot X)$.
			For any $a\in I'$ and $b\in J'$ we have $[a\widecheck{T|_{J'}}]b=aT|_{J'}b=aTb=[a\check{T}]b$, so $\widecheck{T|_{J'}}=\check{T}|_{I'}$.
			Thus the isomorphisms $\Ll_R(J,X\cdot J)\cong\Ll_L(I_J,I_J\cdot X)$ arising from Lemma~\ref{lem-symmetricMultiplierModuleStructure} preserve the inductive limit structures and induce an isomorphism of ${}_\loc X$ and $X_\loc$.
		\end{proof}
	\end{proposition}
	
	Proposition~\ref{prop:leftRightLocAgree} allows us to suppress the left- and right- prefixes of the local multiplier module of a bimodule. 
	In principle these differ as sets, but we shall identify both under the isomorphism specified above.
	This is a problem that does not arise in the case of local multiplier algebras, as these (when considered as Hilbert bimodules) are symmetric, that is, $\Mloc(A)^*=\Mloc(A)$ explicitly. 
	
	\begin{lemma}\label{lem:detectSubbimodsOfLoc}
		Let $X$ be a Hilbert $A$--$B$-bimodule and let $Y\subseteq X_\loc$ be a Hilbert $\Mloc(A)$--$\Mloc(B)$-subbimodule of $X_\loc$. Then $Y=\{0\}$ if and only if $Y\cap X=\{0\}$. 
		\begin{proof}
			Suppose $X\cap Y=\{0\}$.
			We have $Y=X_\loc\cdot s(Y)$ by the classification of Hilbert subbimodules, and so
			\begin{align*}
				\{0\}&=X\cap Y\\
				&=X\cap (X_\loc\cdot s(Y))\\
				&\supseteq (X\cdot(s(Y)\cap B))\cap (X_\loc\cdot s(Y))\\
				&=X\cdot (s(Y)\cap B).
			\end{align*}
			Thus $s(Y)\cap B\subseteq s(X)^\perp$, and so by Lemma~\ref{lem:idealOfMlocContainment} we have $s(Y)\subseteq\Mloc(s(X)^\perp)$.
			Hence $Y=X_\loc\cdot s(Y)\subseteq X_\loc\Mloc(s(X)^\perp)=(X\cdot s(X)^\perp)_\loc=\{0\}$ by Lemma~\ref{lem:locIsMultiplicativeForIdeals}.
		\end{proof}
	\end{lemma}
	
	Lemma~\ref{lem:detectSubbimodsOfLoc} is the bimodule analogue of detection of ideals for an inclusion $A\subseteq\Mloc(A)$. 
	It is worth noting that this argument requires $X$ to be a Hilbert bimodule, since without this condition it is no longer true in general that $Y=X\cdot s(Y)$ (for example, all closed subspaces of a Hilbert space have the same source ideal $\CC$).
	This shall not be a problem for us, as in later sections we shall exclusively examine Hilbert bimodules.
	
	Although the construction of the local multiplier module commutes with right multiplication by an ideal (as in Lemma~\ref{lem:locIsMultiplicativeForIdeals}), this is not true for balanced tensor products of Hilbert bimodules. 
	One particular consequence of this is that if $X$ is a Morita equivalence $A$--$B$-bimodule, it is not always true that $X_\loc$ induces a Morita equivalence between $\Mloc(A)$ and $\Mloc(B)$. 
	This failure occurs because $X_\loc$ need not be full, even when $X$ is.
	\begin{example}\label{eg:localiseHilbSpace}
		Let $X=\Hh$ be an infinite-dimensional Hilbert space. Then $\Hh$ induces a Morita equivalence between $\CC$ and the algebra $\Kk(\Hh)$ of compact operators on $\Hh$, with $s(\Hh)=\CC$ and $r(\Hh)=\Kk(\Hh)$. Clearly, $\Hh_\loc=\Ll(\CC,\Hh)\cong\Hh$ as $\CC$ is simple and unital. Hence $\Hh_\loc\otimes_\CC\Hh_\loc^*\cong\Kk(\Hh)$. However, $(\Hh\otimes_\CC\Hh^*)_\loc\cong\Kk(\Hh)_\loc=\Bb(\Hh)$, the algebra of bounded operators on $\Hh$ (considered as a Hilbert bimodule over $\Mloc(\Kk(\Hh))=\Bb(\Hh)$), which is strictly larger than $\Kk(\Hh)$ since $\Hh$ is infinite-dimensional.
	\end{example}
	
	\begin{proposition}\label{prop:tensorProdLocalMap}
		Let $X$ be a Hilbert $A$--$B$-bimodule and let $Y$ be a Hilbert $B$--$C$-bimodule. There is an isometric bimodule map $\Phi:X_\loc\otimes_{\Mloc(B)}Y_\loc\to(X\otimes_B Y)_\loc$ extending the canonical embedding $X\otimes_B Y\subseteq X_\loc\otimes_{\Mloc(B)}Y_\loc$.
		\begin{proof}
			Fix $J\in\Ii_e(B)$ and $K\in\Ii_e(C)$. 
			Then, there is $J'\in\Ii_e(C)$ such that $J\cdot Y=Y\cdot J'$, and consequently $X\cdot J\otimes_BY\cdot K=X\otimes_B J\cdot Y\cdot K=X\otimes_BY\cdot J'K$. 
			Fix $\xi\in\Ll(J,X\cdot J)$ and $\eta\in\Ll(K,Y\cdot K)$. 
			Denote by $\mult_{J'K,Y}:J'K\otimes_B Y\to J'K\cdot Y$ the unitary map implementing the left multiplication of $B$ on $Y$. 
			Identify $\eta$ with $\eta|_{J'K}\in\Ll(J'K,Y\cdot J'K)$ and define 
			$$\varphi_{\xi,\eta}:=(\xi\otimes 1)\circ\mult_{J'K,Y}^{-1}\circ\eta:J'K\to(X\otimes_BY)J'K.$$
			This is an adjointable map since $\xi\otimes 1, \mult_{J'K,Y},$ and $\eta$ each are. 
			We claim there is a bimodule homomorphism $\Phi:X_\loc\otimes_{\Mloc(B)}Y_\loc\to(X\otimes_B Y)_\loc$ that maps the elementary tensor $\xi\otimes\eta$ to $\varphi_{\xi,\eta}$ for $\xi\in\Ll(J,X\cdot J)$ and $\eta\in\Ll(Y,Y\cdot K)$. 
			First, note from the definition that for $\tau\in M(r(X\cdot J))$, $\omega\in M(J\cap r(Y\cdot J'K))$ and $\sigma\in M(J'K)$, we have $\varphi_{\tau\xi\omega,\eta\sigma}=\tau\varphi_{\xi,\omega\eta}\sigma$. The assignment $(\xi,\eta)\mapsto\varphi_{\xi,\eta}$ thus respects the balanced tensor product structure and $\Phi$ will be a bimodule homomorphism. 
			To see that the assignment is isometric, consider $J'$ as above and fix $\xi_1,\xi_2\in\Ll(J,X\cdot J)$ and $\eta_1,\eta_2\in\Ll(J'K,Y\cdot J'K)$. For $b\in J'K$ we compute
			\begin{align*}
				\langle\varphi_{\xi_1,\eta_1},\varphi_{\xi_2,\eta_2}\rangle(b)&=\varphi_{\xi_1,\eta_1}^*\varphi_{\xi_2,\eta_2}(b)\\
				&=(\eta_1^*\circ\mult_{J'K,Y}\circ\xi_1^*\otimes 1)(\xi_2(a)\otimes y),&\text{where $a\cdot y=\eta_2(b)$,}\\
				&=\eta_1^*(\xi_1^*\xi_2(a)y)\\
				&=\eta_1^*\xi_1^*\xi_2\eta_2(b)\\
				&=\langle\xi_1\otimes\eta_1,\xi_2\otimes\eta_2\rangle(b).
			\end{align*}
			Lastly, for $x\in X$, $y\in Y$, and $c\in C$ we have $\varphi_{x,y}(c)=x\otimes y\cdot c$, so $\Phi$ restricts to the canonical embedding of $X\otimes_BY$ into $(X\otimes_B Y)_\loc$.
		\end{proof}
	\end{proposition}
	
	The map given in Proposition~\ref{prop:tensorProdLocalMap} will in many cases not be adjointable. 
	To see this, recall Example~\ref{eg:localiseHilbSpace} and note that $\Kk(\Hh)$ embeds identically into $\Bb(\Hh)$, but since $\Kk(\Hh)$ is an essential ideal this map cannot have an adjoint. 
	
	The image of $X_\loc\otimes_{\Mloc(B)}Y_\loc$ in the localisation of the tensor product does however `detect' all Hilbert $\Mloc(A)$--$\Mloc(C)$-subbimodules of $(X\otimes_B Y)_\loc$. 
	This follows as $X_\loc\otimes_{\Mloc(B)}Y_\loc$ contains $X\otimes_B Y$, which detects such subbimodules by Lemma~\ref{lem:detectSubbimodsOfLoc}. 
	This then implies that $(X_\loc\otimes_{\Mloc(B)}Y_\loc)^\perp=\{0\}$, giving that $s(X_\loc\otimes_{\Mloc(B)}Y_\loc)$ is an essential ideal of $s((X\otimes_B Y)_\loc)$.
	
	\subsection{Non-triviality conditions for bimodules}
	One may also ask how the construction of the local multiplier module process alters non-triviality conditions for bimodules. 
	We are interested in two particular conditions for Hilbert bimodules: pure outerness and aperiodicity. 
	Pure outerness is of interest to us for the classification of noncommutative Cartan pairs (cf. \cite{E1}, \cite{KM1}). 
	Aperiodicity is explored in \cite{KM2} and grants us two properties of interest: uniqueness of conditional expectations taking value in the injective hull of a $C^*$-algebra, in particular, the local multiplier algebra (\cite[Theorem~1]{F1}, \cite{KM4}), and pure outerness of the local multiplier module of a bimodule (see Proposition~\ref{prop:locPurOutIsAper} ahead). 
	In later sections we shall use this to construct Cartan pairs by applying the local multiplier module construction to a special class of actions that do not initially give rise to Cartan inclusions.
	
	\begin{definition}[{\cite[Definition~4.3]{KM3}}]
		Let $X$ be a Banach $A$-bimodule. We say that $X$ is \emph{purely outer} if the only ideal of $J\triangleleft A$ such that $X\cdot J\cong J$ is $J=\{0\}$.
	\end{definition}
	
	\begin{definition}[{\cite{KM3}}]
		Let $X$ be a normed $A$-bimodule. We say that $x\in X$ satisfies \emph{Kishimoto's condition} if for any $\eps>0$ and any non-zero hereditary subalgebras $D\subseteq A$, there exists $a\in D$ with $a\geq 0$ and $||a||=1$ such that $||axa||<\eps$. We say that $X$ is \emph{aperiodic} if all $x\in X$ satisfy Kishimoto's condition.
	\end{definition}
	
	\begin{lemma}[{\cite[Lemma~5.12]{KM2}}]\label{lem:aperModProps}
		Subbimodules, quotient bimodules, extensions, finite direct sums, and inductive limits of aperiodic normed $A$-bimodules remain aperiodic. If $f:X\to Y$ is a bounded $A$-bimodule homomorphism with dense range and $X$ is aperiodic, then so is $Y$. If $D\subseteq A$ is hereditary, then an aperiodic $A$-bimodule is also aperiodic as a $D$-bimodule. If $J\in\Ii_e(A)$ and $X$ an $A$-bimodule, then $JXJ$ is aperiodic as a $J$-bimodule if and only if $X$ is aperiodic as an $A$-bimodule.
	\end{lemma}

	By \cite[Lemma~5.10]{KM2} no positive non-zero element of $A$ satisfies Kishimoto's condition when considering $A$ as an $A$-bimodule. In particular, no ideals of $A$ are aperiodic other than the zero ideal. Thus, if $X$ is an aperiodic Banach bimodule and $J\triangleleft A$ is an ideal such that $X\cdot J\cong J$, then $J$ is also aperiodic, giving $J=0$.
	This shows that aperiodic Banach bimodules are purely outer. The converse is in general false. However if $X$ is a Hilbert $A$-bimodule and $A$ contains an essential ideal that is simple or of type I, then they are equivalent by \cite[Theorem~8.1]{KM3}.
	
	\begin{proposition}\label{prop:locPurOutIsAper}
		Let $X$ be an aperiodic Hilbert $A$-bimodule. Then $X_\loc$ is a purely outer $\Mloc(A)$-bimodule.
		\begin{proof}
			Let $\tilde{I}\triangleleft \Mloc(A)$ be an ideal such that $X_\loc\cdot\tilde{I}\cong \tilde{I}$. Let $\Phi:X_\loc\cdot\tilde{I}\to\tilde{I}$ be such an isomorphism. Then $\Phi$ restricts to an injective $A$-bimodule map $\phi:X\cap(X_\loc\cdot\tilde{I})\to\Mloc(A)$. 
			The image of $\phi$ is an aperiodic $A$-bimodule by \cite[Lemma~5.12]{KM2}.
			By \cite[Proposition~3.16]{KM4}, the algebra $\Mloc(A)$ contains no non-zero aperiodic $A$-bimodule. 
			Hence the image of $\phi$ is zero, and so $X\cap(X_\loc\cdot\tilde{I})=\{0\}$. Lemma~\ref{lem:detectSubbimodsOfLoc} then gives $\tilde{I}\cong X_\loc\cdot\tilde{I}=\{0\}$.
		\end{proof}
	\end{proposition}
	
	\subsection{A Galois connection for Hilbert modules}
	
	Let $A\subseteq B$ be $C^*$-algebras and let $\Ii(A)$ and $\Ii(B)$ be their respective ideal lattices.
	There are maps $i:\Ii(A)\to\Ii(B)$ and $r:\Ii(B)\to\Ii(A)$ given by
	$$r(J)=A\cap J,\qquad i(I)=BIB,$$
	for $I\in\Ii(A)$ and $J\in\Ii(B)$, where $BIB:=\overline\sspan BIB$ is the ideal generated by $I$ in $B$.
	Green \cite{Green1} observed that these maps form a \emph{monotone Galois connection}.
	That is, the maps $r$ and $i$ have the property that for any $I\in\Ii(A)$ and $J\in\Ii(B)$, the containment $I\subseteq r(J)$ holds if and only if $i(I)\subseteq J$.
	It follows then that the maps $i$ and $r$ also satisfy $i\circ r\circ i=i$ and $r\circ i\circ r=r$.
	The map $i$ preserves joins (sum closure of ideals) and the map $r$ preserves meets (intersections of ideals), and these maps further restrict to mutually inverse isomorphisms
	$$i(\Ii(A))\cong r(\Ii(B)).$$
	
	This argument extends to Hilbert modules.
	Let $X$ be a Hilbert $A$-module and let $Y$ be a Hilbert $B$-module. 
	Suppose there is an isometric embedding $X\subseteq Y$.
	Denote by $\Mod_A(X)$ and $\Mod_B(Y)$ the collections of $A$-submodules of $X$ and $B$-submodules of $Y$ respectively.
	There are analogously defined restriction and induction maps $r:\Mod_B(Y)\to\Mod_A(X)$ and $i:\Mod_A(X)\to\Mod_B(Y)$ given by
	$$r(Z)=Z\cap X,\qquad i(W)=W\cdot B,$$
	where $W\cdot B=\overline\sspan W\cdot B$ is the Hilbert $B$-module generated by $W$ in $Y$.
	One then readily checks
	\begin{align*}
		W\subseteq r(Z)&\implies i(W)=W\cdot B\subseteq r(Z)\cdot B\subseteq Z,\text{ and}\\
		i(W)\subseteq Z&\implies i(W)\cap X\subseteq Z\cap X\\
		&\implies W\subseteq r(W),
	\end{align*}
	so these maps also induce a monotone Galois correspondence.
	A similar argument works for inclusions $A_i\subseteq B_i$ for $i=1,2$ and Hilbert bimodules $X\subseteq Y$, where $X$ is a Hilbert $A_1$--$A_2$-bimodule and $Y$ is a Hilbert $B_1$--$B_2$-bimodule.
	Thus we have shown that $r$ is a right adjoint to $i_-$.\\
	
	In the setting of the inclusion $A\subseteq \Mloc(A)$, have an alternative way to induce ideals of $\Mloc(A)$ and Hilbert $\Mloc(A)$-(bi)modules.
	If $I\triangleleft A$ is an ideal and if $X$ is a Hilbert $A$-module, we gain an ideal $i_+(I):=\Mloc(I)$ of $\Mloc(A)$ and may also induce the Hilbert $\Mloc(A)$-(bi)module $i_+(X):=X_\loc$ from $X$.
	These maps $i_+$ do not however induce a Galois correspondence with restriction $r$.
	One counterexample is to consider the compact operators $\Kk(\Hh)$ as in Example~\ref{eg:localiseHilbSpace}.
	Since the local multiplier algebra of $\Kk(\Hh)$ is the bounded operators $\Bb(\Hh)$, and $\Kk(\Hh)$ is also an ideal in $\Bb(\Hh)$, we see that $\Kk(\Hh)=\Kk(\Hh)\cap\Bb(\Hh)=r(\Kk(\Hh))$ but $i_+(\Kk(\Hh))=\Bb(\Hh)$ is not contained in $\Kk(\Hh)$ if $\Hh$ has infinite dimension, violating the Galois correspondence condition.
	
	If we restrict to the Hilbert subbimodules $W\subseteq X$, $Z\subseteq Y$ satisfying $W=r(i_+(W))$ and $Z=i_+(r(Z))$, we do however gain
	$$r(Z)\subseteq W\iff Z\subseteq i_+(W).$$
	Restricting to this class of ideals, we gain a Galois correspondence and so $i_+$ is then the right adjoint to $r$ for this restricted class.
	
	\section{Inducing inverse semigroup actions on the local multiplier algebra}
	
	We shall give the definition of an inverse semigroup action on a $C^*$-algebra $A$ by Hilbert bimodules as defined and used in \cite{BM1}, \cite{KM1}, \cite{KM2}, and then apply our local multiplier module construction.
	This gives rise to a Fell bundel over the acting inverse semigroup, which after some refinement, gives an action of an inverse semigroup on $\Mloc(A)$.
	
	\subsection{Inverse semigroup actions and crossed products}
	
	Throughout the rest of this article, $A$ shall denote a $C^*$-algebra and $S$ shall denote a unital inverse semigroup. 
	There is a canonical partial order on $S$ given by $t\leq u$ if and only if $t=ut^*t$ for $u,t\in S$.
	
	\begin{definition}[cf.\ {\cite[Definition~4.7]{BM1}}]
		An \emph{action} $\Ee=(\Ee_t,\mu_{t,u})_{t,u\in S}$ of $S$ on $A$ by Hilbert bimodules consists of
		\begin{itemize}
			\item Hilbert $A$-bimodules $\Ee_t$ for each $t\in S$; and
			\item bimodule isomorphisms $\mu_{t,u}:\Ee_t\otimes_A\Ee_u\to\Ee_{tu}$;
		\end{itemize}
		satisfying
		\begin{enumerate}[(i)]
			\item $\Ee_1=A$, with the canonical $A$-bimodule structure;
			\item the maps $\mu_{1,t}:A\otimes_A\Ee_t\to\Ee_t$ and $\mu_{t,1}:\Ee_t\otimes_A A\to\Ee_t$ are the canonical isomorphisms coming from the respective left and right actions of $A$ on $\Ee_t$; and
			\item associativity: for all $t,u,v\in S$, the following diagram commutes:
			\[
   				\begin{tikzpicture}[baseline=(current bounding box.west)]
			    	\node (1) at (0,1) {$(\Ee_t\otimes_A \Ee_u) \otimes_A \Ee_v$};
      				\node (1a) at (0,0) {$\Ee_t\otimes_A (\Ee_u \otimes_A \Ee_v)$};
      				\node (2) at (5,1) {$\Ee_{tu} \otimes_A \Ee_v$};
      				\node (3) at (5,0) {$\Ee_t\otimes_A \Ee_{uv}$};
      				\node (4) at (7,.5) {$\Ee_{tuv}$};
      				\draw[<->] (1) -- node[swap] {ass.} (1a);
      				\draw[cdar] (1) -- node {$\mu_{t,u}\otimes_A \Id_{\Ee_v}$} (2);
      				\draw[cdar] (1a) -- node[swap] {$\Id_{\Ee_t}\otimes_A\mu_{u,v}$} (3);
     				\draw[cdar] (3.east) -- node[swap] {$\mu_{t,uv}$} (4);
      				\draw[cdar] (2.east) -- node {$\mu_{tu,v}$} (4);
    			\end{tikzpicture}
    		\]
		\end{enumerate}
	\end{definition}	
	
	We shall now give the necessary concepts to build the full and essential crossed products for inverse semigroup actions. These are the same definitions as in \cite[Section~2.2]{KM1} and \cite[Section~4]{KM2}, and we refer the reader there for a more in depth explanation.
	
	If $t\leq u$ for $t,u\in S$ and $\Ee$ is an action of $S$ on $A$, then there is an inclusion map $\Ee_t\hookrightarrow\Ee_u$ gained from the multiplication maps. These inclusions then restrict to isomorphisms $j_{u,t}:\Ee_t\to\Ee_{u}\cdot s(\Ee_t)=r(\Ee_t)\cdot\Ee_u$. For each $v\leq t,u$, the maps $j_{t,v}$ and $j_{u,v}$ induce an isomorphism $\vartheta_{t,u}^v:\Ee_u\cdot s(\Ee_v)\to\Ee_t\cdot s(\Ee_v)$ by defining $\vartheta_{t,u}^v:=j_{t,v}\circ j_{u,v}^{-1}$. Define
	$$I_{t,u}=\overline{\sum_{v\leq t,u}s(\Ee_v)},$$
	the closed ideal generated by $s(\Ee_v)$ for all $v\leq t,u$. We call the ideal $I_{t,u}$ the \emph{intersection ideal for $t,u$}. 
	This is contained in $s(\Ee_u)\cap s(\Ee_t)$ and the inclusion may be strict. 
	There is a unique Hilbert bimodule isomorphism $\vartheta_{t,u}:\Ee_u\cdot I_{t,u}\to\Ee_v\cdot I_{t,u}$ which for each $v\leq t,u$ restricts to $\vartheta_{t,u}^v$ on $\Ee_t\cdot s(\Ee_v)$ by \cite[Lemma~2.4]{BM1}. 
	The \emph{algebraic crossed product} $A\rtimes_\alg S$ is defined as the quotient of $\bigoplus_{t\in S}\Ee_t$ by the linear span of $\vartheta_{u,t}(\xi)\delta_u-\xi\delta_t$ for $t,u\in S$ and $\xi\in\Ee_t\cdot I_{t,u}$. 
	There is a ${}^*$-algebra structure on $A\rtimes_\alg S$ with multiplication and involution induced by the maps $\mu_{t,u}$ and the involutions $\Ee_t^*\to\Ee_{t^*}$. 
	There is then a maximal $C^*$-norm on $A\rtimes_\alg S$.	
	
	\begin{definition}[{\cite[Definition~2.7]{KM1}}]
		The \emph{full crossed product} $A\rtimes S$ of the action $\Ee$ is defined as the maximal $C^*$-completion of the ${}^*$-algebra $A\rtimes_\alg S$.
	\end{definition}
	
	We also define some non-triviality conditions for inverse semigroup actions.
	Kwa\'sniewski and Meyer (\cite{KM1}, \cite{KM2}) show that actions satisfying these conditions give rise to inclusions of $C^*$-algebras satisfying certain maximality conditions that are useful in the analysis of Cartan-like pairs.
	
	\begin{definition}[{\cite[Definitions~6.1, 6.9]{KM2}}]
		Let $\Ee$ be an action of an inverse semigroup $S$ on $A$. We say that the action is \emph{purely outer} if $\Ee_t\cdot I_{1,t}^\perp$ is a purely outer $A$-bimodule for each $t\in S$. We say that the action is \emph{aperiodic} if $\Ee_t\cdot I_{1,t}^\perp$ is an aperiodic $A$-bimodule for each $t\in S$.
	\end{definition}
	
	A definition of the reduced crossed product for such an action can be found in \cite[Section~2.2]{KM1}. This involves the construction of a weak conditional expectation $E:A\rtimes S\to A''$ taking values in the enveloping von Neumann algebra of $A$. One may also consider \emph{essentially defined conditional expectations} or \emph{local expectations} that instead take values in the local multiplier algebra $\Mloc(A)$ of $A$.
	
	Recall that for a $C^*$-algebra $A$ we denote Hamana's injective hull by $I(A)$ \cite{H1}.
	
	\begin{definition}[{\cite[Definitions~3.1,3.7]{KM2}}]
		Let $A\subseteq B$ be a $C^*$-inclusion.
		A \emph{generalised expectation} for $A\subseteq B$ consists of another $C^*$-inclusion $A\subseteq\tilde{A}$ and a completely positive contractive linear map $E:B\to\tilde{A}$ such that $E$ restricts to the identity on $A$. 
		We say that $E$ is \emph{faithful} if $E(b^*b)=0$ implies $b=0$ for all $b\in B$, and $E$ is \emph{almost faithful} if $E((ba)^*ba)=0$ for all $a\in B$ implies $b=0$ for all $b\in B$.
		If $\tilde{A}=\Mloc(A)$ we call $E$ an \emph{essentially defined conditional expectation} or a \emph{local expectation}.
		If $\tilde{A}=I(A)$ then we call $E$ a \emph{pseudo-expectation}.
	\end{definition}
	
	Since $\Mloc(A)$ canonically embeds in $I(A)$ by a result of Frank \cite{F1} we can consider $\Mloc(A)\subseteq I(A)$.
	With this one may consider local expectations as a special case of pseudo-expectations.
	
	The dense subalgebra $A\rtimes_\alg S$ of $A\rtimes S$ is exactly the span of the bimodules $\Ee_t$ under the canonical inclusions $\Ee_t\to A\rtimes_\alg S$, $\xi\mapsto[\xi\delta_t]$. 
	Each of these maps is injective, so we identify each $\Ee_t$ with its image in $A\rtimes_\alg S$. 
	We define a local expectation $EL:A\rtimes S\to\Mloc(A)$ as follows
	Let $t\in S$. 
	Recalling the isomorphism $\vartheta_{1,t}:\Ee_t\cdot I_{1,t}\to I_{1,t}$, each $\xi\in\Ee_t$ defines an element of $\Mm(I_{1,t})\subseteq\Mloc(A)$ via $EL(\xi):a\mapsto\vartheta_{1,t}(\xi\cdot a)$. 
	By \cite[Proposition~4.4]{KM2} this extends to a local expectation $EL:A\rtimes S\to\Mloc(A)$, which we call the \emph{canonical local expectation}. 
	We denote by $\Nn_{EL}$ the largest ideal contained in $\ker(EL)$.
	
	\begin{definition}
		The \emph{essential crossed product} is defined as the quotient 
		$$A\rtimes_\ess S:=(A\rtimes S)/\Nn_{EL}.$$ 
		The local expectation $EL$ descends to a local expectation $A\rtimes_\ess S\to\Mloc(A)$, which we also denote by $EL$.
	\end{definition}
	
	\begin{theorem}[{\cite[Theorem~4.12]{KM2}}]
		The canonical local expectation $EL:A\rtimes_\ess S\to\Mloc(A)$ is faithful.
	\end{theorem}
	
	Renault \cite{R1}, and earlier Kumjian \cite{K1} considered inclusions of $C^*$-algebras $A\subseteq B$ such that the larger algebra $B$ is densely spanned by elements that normalise $A$ in the following sense: $n\in B$ is a normaliser if $n^*An,nAn^*\subseteq A$. This has since been taken as a standard assumption for Cartan pair-like objects in \cite{E1}, \cite{KM1}, and many others. For an action $\Ee$ of a unital inverse semigroup $S$ on a $C^*$-algebra $A$ (or more generally a Fell bundle over $S$ with unit fibre $A$), the inclusion $A\subseteq A\rtimes S$ has this property. This is because each $\Ee_t$ carries an $A$-bimodule structure, and together they span a dense subspace of $A\rtimes S$. 
	
	\begin{definition}
		Let $A\subseteq B$ be an inclusion of $C^*$-algebras. We call the inclusion \emph{regular} if the set of normalisers $N(A,B)=\{n\in B: n^*An,nAn^*\subseteq A\}$ spans a dense subspace of $B$. Closed subspaces $M\subseteq N(A,B)$ such that $AM,MA\subseteq M$ are called \emph{slices}, and the collection of slices for the inclusion $A\subseteq B$ is denoted $\Ss(A,B)$. A \emph{subslice} $N$ of $M\in\Ss(A,B)$ is a slice $N$ contained in $M$.
	\end{definition}
	
	If $A\subseteq B$ is a non-degenerate inclusion then $M^*M, MM^*\subseteq A$ for any slice $M\in\Ss(A,B)$, giving each slice a Hilbert $A$-bimodule structure with inner products induced from the multiplication in $B$. The set $\Ss(A,B)$ then becomes an inverse semigroup with operation $M\cdot N:=\cspan MN$ and ${}^*$-operation given by the adjoint in $B$. For brevity of notation we write $MN$ to denote the closed linear span of products of elements in $M$ and $N$. This gives rise to the tautological action of $\Ss(A,B)$ on $A$, where the bimodules for the action are the slices, and the multiplication isomorphisms are induced by the multiplication in $B$. In the case where one has a closed and purely outer action $\Ee$ of $S$ on $A$, Kwa\'sniewski and Meyer \cite{KM1} showed one can reconstruct slices for the inclusion $A\subseteq A\rtimes_\red S$ from the bimodules $\Ee_t$.
	
	In general if $A\subseteq B$ is a regular non-degenerate inclusion then the slice inverse semigroup $\Ss(A,B)$ acts	tautologically on $A$ via the action $\Ee_X=X$ for slices $X\in\Ss(A,B)$, and multiplication maps given by the multiplication in $B$.
	
	We then gain the following classification of certain $C^*$-inclusions in terms of actions by inverse semigroups.
	
	\begin{theorem}\label{thm:mainAperIso}
		Let $A\subseteq B$ be a non-degenerate inclusion of $C^*$-algebras, and let $E:B\to I(A)$ be a faithful pseudo-expectation. 
		Suppose there is a densely spanning inverse subsemigroup $S\subseteq\Ss(A,B)$ that acts aperiodically on $A$. 
		Then there is an isomorphism $\varphi:A\rtimes_\ess S\to B$ that restricts to the identity on $A$ and entwines $E$ with the canonical local expectation $EL$. 
		In particular, the expectation $E$ takes values in $\Mloc(A)$ and the inclusion $A\subseteq B$ is aperiodic.
		\begin{proof}
			The inclusion $A\subseteq A\rtimes S$ is aperiodic by \cite[Proposition~6.3]{KM2}, so the local expectation $EL:A\rtimes S\to \Mloc(A)\subseteq I(A)$ is the unique pseudo-expectation for the inclusion by \cite[Theorem~3.6]{KM4}.
			Let $\Phi:A\rtimes S\to B$ be the canonical ${}^*$-homomorphism.
			This is surjective as it spans each of the slices in $S$.
			Then $E\circ\Phi:A\rtimes S\to I(A)$ is a pseudo-expectation, and so is equal to $EL$ by uniqueness.
			
			The map $\Phi$ then descends to an isomorphism $\varphi:A\rtimes S/\ker(\Phi)\to B$, so it suffices to show that the kernel of $\Phi$ is $\Nn_{EL}$, the largest ideal contained in the kernel of $EL$, as taking the quotient by this ideal gives the essential crossed product.
			Fix $x\in\ker(\Phi)$. 
			Then $EL(x^*x)=E(\Phi(x^*x))=0$ so $x\in\Nn_{EL}$.
			Conversely, $x\in\Nn_{EL}$ if and only if $0=EL((xy)^*xy)$ by \cite[Proposition~3.5]{KM2}, so $0=EL((xy)^*xy)=E(\Phi(xy)^*\Phi(xy))$.
			Since $E$ is faithful, it is almost faithful by \cite[Corollary~3.7]{KM2}, and since $\Phi$ is surjective we see that $\Phi(x)=0$.
			Thus $\Phi$ descends to an isomorphism $\varphi:A\rtimes_\ess S\to B$.
			The inclusion $A\subseteq A\rtimes_\ess S$ is then aperiodic since the quotient map $q:A\rtimes S\to A\rtimes_\ess S$ descends to a bounded surjective bimodule map $(A\rtimes S)/A\to (A\rtimes_\ess S)/A$, and the image of an aperiodic bimodule is aperiodic by Lemma~\ref{lem:aperModProps}.
		\end{proof}
	\end{theorem}	
	
	The universal property of $I(A)$ ensures that for any $C^*$-inclusion $A\subseteq B$ there is always a pseudo-expectation $E:B\to I(A)$ extending the identity map on $A$.
	Under the conditions of Theorem~\ref{thm:mainAperIso} we then gain that $E$ must in fact take values in $\Mloc(A)$, as the canonical pseudo-expectation associated to the essential crossed product does.
	
	\subsection{The dual groupoid to an inverse semigroup action}
	
	If $\Ee$ is an action of $S$ on $A$, there is an induced action $\hat\Ee=(\hat\Ee_t)_{t\in S}$ of $S$ on $\hat{A}$ such that $\hat\Ee_t:\widehat{s(\Ee_t)}\to\widehat{r(\Ee_t)}$ for each $t\in S$. The construction of the transformation groupoid $\hat{A}\rtimes S$ associated to $\Ee$ can be found in \cite[Section~2.3]{KM1}.
	
	\begin{definition}[{\cite{KM2}}]
		We call $\hat\Ee$ the \emph{dual action} to the action $\Ee$ of $S$ on $A$. The transformation groupoid $\hat{A}\rtimes S$ is called the \emph{dual groupoid}.
	\end{definition}
	
	The unit space of the dual groupoid $\hat{A}\rtimes S$ is homeomorphic to $\hat{A}$ via the map $\hat{A}\to(\hat{A}\rtimes S)^{(0)}$, $[\pi]\mapsto[1,[\pi]]$. 
	We often identify $\hat{A}$ and $(\hat{A}\rtimes S)^{(0)}$ under this map.
	By \cite[Proposition~2.17]{KM1}, an action $\Ee:S\curvearrowright A$ is closed if and only if the unit space $\hat{A}$ of $\hat{A}\rtimes S$ is a closed subset.

	By construction of the essential crossed product, there are injective $A$-bimodule maps $\Ee_t\to A\rtimes_\ess S$ where the image of each $\Ee_t$ is a slice for this inclusion.
	Let $B:= A\rtimes_\ess S$.
	We gain a canonical homomorphism $S\to\Ss(A,B)$ by mapping $t\in S$ to the corresponding image of $\Ee_t$ in $B$.
	We then gain a canonical groupoid homomorphism $\phi:\hat{A}\rtimes S\to\hat{A}\rtimes\Ss(A,B)$ via $\phi[t,[\pi]]=[\Ee_t,[\pi]]$.
	The inverse semigroup $S$ also canonically maps into the bisection inverse semigroup $\Bis(\hat{A}\rtimes\Ss(A,B))$ by mapping $t\in S$ to the bisection $\{[\Ee_t,[\pi]]:[\pi]\in \widehat{s(\Ee_t)}\}$.
	In general these homomorphisms are not injective.
	For example, if $A_1$ and $A_2$ are $C^*$-algebras, consider the action generated by a self-inverse automorphism on $A_1$, and extend this to an automorphism $\alpha$ on $A:=A_1\oplus A_2$ by acting trivially on $A_2$.
	Then the germ relation for a group action is trivial since groups only have one idempotent: the unit.
	Thus over every irreducible representation of $A_2$ there is non-trivial isotropy, but if we consider $\hat{A}\rtimes\Ss(A,A\rtimes_\ess \ZZ_2)$, we see that $[\alpha,[\pi]]=[\Id,[\pi]]$ for all $[\pi]\in A_2$ as $\alpha|_C=\Id|_C$.
	These two germs differ in $\hat{A}\rtimes\ZZ_2$, since the only germ relations we may take there are trivial, so $[\alpha,[\pi]]\neq[\Id,[\pi]]$ in $\hat{A}\rtimes\ZZ_2$ for all $[\pi]\in\hat{A}$.
	
	This can be remedied by assuming that our action has enough idempotents.
	
	\begin{lemma}\label{lem-enoughIdemsGivesInjGrpdMap}
		Let $\Ee$ be a purely outer action of $S$ on $A$ such that for each the open sets $\hat\Ee_e\subseteq\hat{A}$ for idempotents $e\in S$ form a basis for the topology on $\hat{A}$.
		Write $B:= A\rtimes_\ess S$.
		Then the canonical homomorphism $\phi:\hat{A}\rtimes S\to\hat{A}\rtimes\Ss(A,B)$ is injective.
		\begin{proof}
			Fix $[t,[\pi]]\in \hat{A}\rtimes S$ with $[\hat\Ee_t,[\pi]]=[A,[\pi]]$, that is, $\phi[t,[\pi]]$ is a unit in $\hat{A}\rtimes\Ss(A,B)$.
			Then there is an open subset $U\subseteq\hat{A}$ such that $\hat\Ee_t|_U=\Id_U$.
			Let $e\in S$ be an idempotent such that $\hat\Ee_e\subseteq U$ and $[\pi]\in\Ee_e$.
			Then $\Ee_t\cdot \Ee_e=\Ee_e$ giving  $\Ee_e=\Ee_t\cdot s(\Ee_e)$, and we see that $e=te$ since the action $\Ee$ is purely outer.
			Thus $[t,[\pi]]=[e,[\pi]]$ is a unit since $e$ is an idempotent.
		\end{proof}
	\end{lemma}
	
	We note that by \cite[Theorem~7.2]{BE1} we can without loss of generality assume our inverse semigroup has enough idempotents, so that the criteria of Lemma~\ref{lem-enoughIdemsGivesInjGrpdMap} are always satisfied.
	
	\subsection{Extending actions to local multiplier algebras}
	
	The results of Exel in \cite{E1} and Kwa\'sniewski--Meyer in \cite{KM1} require that the action of $S$ on $A$ is \emph{closed}, that is, the conditional expectation $EL$ takes values in $A$.
	Since for our purposes the canonical expectation $EL$ need not take values in $A$, we are unable to immediately use these results.
	To circumvent this, we construct a closed action on $\Mloc(A)$ from the action on $A$, then show that the inclusion $A\subseteq A\rtimes_\ess S$ embeds into the crossed product associated to this extended action.
	
	Unfortunately, one cannot simply take an action $\Ee$ of $S$ on $A$, replace the modules with the local multiplier counterparts, and then gain an action on $\Mloc(A)$. 
	The obstruction to this is that the map in Proposition~\ref{prop:tensorProdLocalMap} is not always an isomorphism, so we do not always have $(\Ee_t)_\loc\otimes_{\Mloc(A)}(\Ee_u)_\loc\cong(\Ee_t\otimes_A\Ee_u)_\loc$.
	Our solution is to instead create a non-saturated Fell bundle with these local multiplier modules, and then gain a saturated Fell bundle using \cite[Theorem~7.2]{BE2}, which is then equivalent to an inverse semigroup action on $\Mloc(A)$.
	
	\begin{definition}[{\cite[Definition~2.10]{BE1}}]\label{defn:fellBundleOverInvSemGrp}
		Let $S$ be an inverse semigroup. A \emph{Fell bundle} over $S$ is a collection $\Aa=(\Aa_t)_{t\in S}$ of Banach spaces $\Aa_t$ together with a multiplication $\cdot:\Aa_t\times\Aa_u\to\Aa_{tu}$ for all $t,u
		\in S$, an involution ${}^*:\Aa_t\to\Aa_{t^*}$ for each $t\in S$, linear maps $j_{t,u}:\Aa_u\to\Aa_t$ for $t,u
		\in S$ with $u\leq t$ satisfying the following:
		\begin{enumerate}[(i)]
			\item the multiplication is bilinear $\Aa_t\times\Aa_u\to\Aa_{tu}$ for all $t,u\in S$;
			\item the multiplication is associative;
			\item $||a\cdot b||\leq ||a||\cdot ||b||$ for all $a,b\in\bigcup_{t\in S}\Aa_t$;
			\item ${}^*$ is conjugate linear on each $\Aa_t$;
			\item $(a^*)^*=a$, $||a^*||=||a||$, and $(a\cdot b)^*=b^*\cdot a^*$ for all $a,b\in\bigcup_{t\in S}\Aa_t$;
			\item $||a^*a||=||a||^2$ and $a^*a$ is a positive element in the $C^*$-algebra $\Aa_{t^*t}$ for all $t\in S$ and $a\in\Aa_t$;
			\item $j_{t,u}$ is an isometric linear map for all $t,u\in S$ with $u\leq t$;
			\item if $v\leq u\leq t$ then $j_{t,v}=j_{t,u}\circ j_{u,v}$;
			\item if $s\leq t$ and $u\leq v$ in $S$ then $j_{t,s}(a)\cdot j_{v,u}(b)=j_{tv,su}(a\cdot b)$ for all $a\in\Aa_s$ and $b\in \Aa_v$;
			\item if $s\leq t$ then $j_{t,s}(a)^*=j_{t^*,s^*}(a^*)$ for all $a\in\Aa_s$. 
		\end{enumerate} 
		If $\Aa_s\cdot\Aa_t$ spans a dense subspace of $\Aa_{st}$ for all $s,t\in S$, we say that the Fell bundle $\Aa$ is \emph{saturated}. If $S$ is unital we call $\Aa_1$ the \emph{unit fibre} of the Fell bundle $\Aa$.
	\end{definition}
		
	One can build a $C^*$-algebra out of the sections of a Fell bundle.	
	
	\begin{definition}[{\cite[Definition~3.4]{E1}, \cite[Definition~2.7]{KM1}}]
		Let $\Aa$ be a Fell bundle over a unital inverse semigroup $S$ with unit fibre $\Aa_1=A$. Let $\Ll(\Aa):=\bigoplus_{t\in S}\Aa_t$ and $N:=\sspan\{a_s\delta_s-j_{t,s}(a_s)\delta_t:s,t\in S$, $s\leq t$, $a_s\in\Aa_s\}$. The \emph{full cross sectional $C^*$-algebra} is $C^*(\Aa)$ is defined as the maximal $C^*$-completion of $\Ll(\Aa)/N$ with multiplication and involution inherited from the Fell bundle.
	\end{definition}	
	
	Each fibre $\Aa_t$ of the Fell bundle $\Aa$ embeds canonically in $\Ll(\Aa)/N$ via $\Aa_t\ni x\mapsto [x\delta_t]\in\Ll(\Aa)/N$, and so embeds in the full $C^*$-algebra $C^*(\Aa)$. Identifying each fibre with its image in $C^*(\Aa)$ we see that $\Ee_t\cdot I_{t,u}=\Ee_u\cdot I_{t,u}=\Ee_t\cap\Ee_u$ for each $t,u\in S$. In particular, we have $I_{1,t}=\Ee_t\cap A$.
	
	\begin{proposition}\label{prop:localisedFellBundle}
		Let $A$ be a $C^*$-algebra, $S$ an inverse semigroup, and $\Ee=(\Ee_t,\mu_{t,u})_{t,u\in S}$ an action of $S$ on $A$ by Hilbert bimodules. There exists a Fell bundle $\Aa$ over $S$ with unit fibre $\Mloc(A)$ such that $\Aa_{t}=(\Ee_t)_\loc$ for each $t\in S$.
		\begin{proof}
			Define $\Aa_t:=(\Ee_t)_\loc$. For $t,u\in S$ we define the multiplication by
			$$\xi_t\cdot\eta_u:=\tilde\mu_{t,u}(\xi_t\otimes_{\Mloc(A)}\eta_u),\qquad \xi_t\in\Aa_t,\eta_u\in\Aa_u,$$
			where $\tilde\mu_{t,u}$ is the map induced by $\mu_{t,u}$ using Proposition~\ref{prop:tensorProdLocalMap}.
			By Lemma~\ref{lem:starLocCommute} and Proposition~\ref{prop:leftRightLocAgree} we canonically identify $(\Ee_t)_\loc^*$ and $(\Ee_{t^*})_\loc$, and can define the ${}^*:\Aa_t\to\Aa_{t^*}$ by mapping each $\xi_t\in(\Ee_t)_\loc$ to its corresponding adjoint $\xi_t^*\in(\Ee_t)_\loc^*=(\Ee_{t^*})_\loc$.
			
			One readily checks that this structure satisfies the axioms of a Fell bundle.
		\end{proof}
	\end{proposition}

	\begin{definition}
		We call the Fell bundle $\Aa$ defined in Proposition~\ref{prop:localisedFellBundle} the \emph{induced local Fell bundle} of the action $\Ee$.
	\end{definition}
	
	The main use of the induced local Fell bundle will be that it `closes' the action.
	By this, we mean the associated conditional expectation for this Fell bundle shall be a genuine conditional expectation $C^*(\Aa)\to\Mloc(A)$ rather than a generalised local expectation taking values in $\Mloc(\Mloc(A))$ (note here that the term `genuine' is in reference to the inclusion $\Mloc(A)\subseteq C^*(\Aa)$).
	This is however not immediately true, but it can be done without loss of generality.
	To show this, we recall that any such inverse semigroup action can be refined to one with ``enough idempotents''.
	
	\begin{proposition}[{\cite[Proposition~5.2]{KM1}}]\label{prop-refinedAction}
		Let $\Ee$ be an action of a unital inverse semigroup $S$ on a $C^*$-algebra $A$.
		Let $\iota:S\to\Bis(\hat{A}\rtimes S)$ be the canonical inclusion $\iota(t)=U_t:=\{[\hat\Ee_t,[\pi]]:[\pi]\in\widehat{s(\Ee_t)}\}$.
		There is an action $\bar\Ee$ of $\bar{S}:=\Bis(\hat{A}\rtimes S)$ on $A$ such that $\bar\Ee_{\iota(t)}=\Ee_t$ for all $t\in S$ and a canonical isomorphism $A\rtimes S\cong A\rtimes\Bis(\hat{A}\rtimes S)$ that restricts to the identity on each $\Ee_t=\bar\Ee_{\iota(t)}$.
		If the map $\iota$ is surjective, we say that $\Ee$ is \emph{fine}.
	\end{proposition}
	
	The action $\bar\Ee$ in Proposition~\ref{prop-refinedAction} is called the \emph{refinement} of $\Ee$.
	The proof of Proposition~\ref{prop-refinedAction} in \cite{KM1} shows that each bimodule $\bar\Ee_U$ for $U\in\Bis(\hat{A}\rtimes S)$ is of the form
	$$\bar\Ee_U=\overline{\sum_{i\in C} \Ee_{t_i}J_i},$$
	for some $t_i\in S$ and $J_i\triangleleft A$, where $C$ is some indexing set.
	
	\begin{lemma}\label{lem-refinedActionAperiodicity}
		The action $\Ee$ is aperiodic if and only if the refinement $\bar\Ee$ is aperiodic.
		\begin{proof}
			If the refined action $\bar\Ee$ is aperiodic then $\Ee_t\cdot I_{1,t}^\perp=\bar\Ee_{U_t}\cdot I_{1,U_t}^\perp$ is aperiodic for all $t\in S$.
			Thus $\Ee$ is aperiodic.
			
			Suppose the action $\Ee$ is aperiodic.
			Fix $U\in\Bis(\hat{A}\rtimes S)$ and write $\bar\Ee_U=\overline{\sum_{i\in C}\Ee_{t_i}\cdot J_i}$ for some $t_i\in S$ and $J_i\triangleleft A$, where $i\in C$ and $C$ is some indexing set.
			For each $U\in\Bis(\hat{A}\rtimes S)$ the supremum over $V\leq 1,U$ in $\Bis(\hat{A}\rtimes S)$ is the bisection $U\cap\hat{A}$. 
			Thus the intersection ideal $I_{1,U}$ is given by $\bar\Ee_{U\cap\hat{A}}=\overline{\sum_i \Ee_{t_i}\cdot J_i}\cap A$.
			For each $i\in C$ we have
			\begin{align*}
				J_i\cdot I_{1,U}^\perp &=J_i\cdot\left(\overline{\sum_{k\in C}\Ee_{t_k}\cdot J_k}\cap A\right)^\perp\\
				&\subseteq J_i\cdot (\Ee_{t_i}\cdot J_i\cap A)^\perp\\
				&=J_i\cdot (\Ee_{t_i}\cap J_i)^\perp.
			\end{align*}
			Thus for each $i\in C$ we have
			\begin{align*}
				\Ee_{t_i}\cdot J_i\cdot I_{1,U}^\perp &\subseteq\Ee_{t_i}\cdot J_i\cdot (\Ee_{t_i}\cap J_i)^\perp\\
				&=\Ee_{t_i}\cdot J_i\cdot I_{1,t}^\perp\\
				&\subseteq\Ee_{t_i}\cdot I_{1,t_i}^\perp.
			\end{align*}
			The bidmodule $\Ee_{t_i}\cdot I_{1,t_i}^\perp$ is aperiodic as the action $\Ee$ is aperiodic, and so we see that each $\Ee_{t_i}\cdot J_i\cdot I_{1,U}^\perp$ is an aperiodic bimodule.
			The span of these bimodules over $i\in C$ is a dense subspace of $\bar\Ee_U\cdot I_{1,U}^\perp$.
			Thus the bimodule $\bar\Ee_U\cdot I_{1,U}^\perp$ is aperiodic by \cite[Lemma~4.2]{KM3}, and so the action $\bar\Ee$ is aperiodic.
		\end{proof}
	\end{lemma}
	
	\begin{corollary}\label{cor-refinedEssCrossProd}
		The canonical isomorphism $A\rtimes_\Ee S\cong A\rtimes_{\bar\Ee}\bar{S}$ entwines the canonical local expectations associated to each action, and hence the isomorphism descends to an isomorphism $A\rtimes_{\Ee,\ess}S\cong A\rtimes_{\bar\Ee,\ess}\bar{S}$ of the essential crossed products.
		\begin{proof}
			The canonical local conditional expectations agree on each $\Ee_t=\bar\Ee_{U_t}$ for each $t\in S$, and these span dense subspaces of both $A\rtimes_\Ee S$ and $A\rtimes_{\bar\Ee}\bar{S}$.
			Thus the isomorphism descends to the essential crossed products.
		\end{proof}
	\end{corollary}
	
	Proposition~\ref{prop-refinedAction}, Lemma~\ref{lem-refinedActionAperiodicity}, and Corollary~\ref{cor-refinedEssCrossProd} together allow us to, without loss of generality, assume that any action we wish to consider is refined.
	In particular, we gain the intersection property that allows us to write any intersection ideal associated to $t,u\in S$ as $I_{t,u}=\Ee_{t\cap u}$ for some $t\cap u\in S$.
	From this point on we identify an action $\Ee$ with its image in the refinement in this way, or we may assume the action $\Ee$ is fine.

	\begin{lemma}\label{lem:intersectIdealLocalised}
		Let $\Ee:S\curvearrowright A$ be a fine action, so that $S\cong\Bis(\hat{A}\rtimes S)$.
		For each $t\in S$ let $\Ii_{1,t}:=\overline{\sum_{v\leq 1,t}s(\Aa_v)}$ be the intersection ideal for $1,t\in S$ for the induced local Fell bundle $\Aa$. 
		Assume that the action $\Ee$ is fine.
		Then $\Ii_{1,t}=\Mloc(I_{1,t})$.
		\begin{proof}
			Identifying $S$ with $\Bis(\hat{A}\rtimes S)$ we see that $\Ii_{1,t}=\tilde\Ee_{t\cap\hat{A}}=(\Ee_{t\cap\hat{A}})_\loc=(I_{1,t})_\loc=\Mloc(I_{1,t})$.
		\end{proof}
	\end{lemma}
	
	\begin{corollary}\label{cor:condExpGen}
		The canonical weak conditional expectation $\tilde{E}:C^*(\Aa)\to\Mloc(A)''$ and local expectation $\tilde{EL}:C^*(\Aa)\to\Mloc(\Mloc(A))$ are genuine \textup{(}that is, take values in $\Mloc(A)$, and agree\textup{)}.
		\begin{proof}
			Using the characterisation of $\tilde{EL}$ in \cite[Lemma~4.5]{BEM1} and the definition of $\tilde{EL}$, we see that for $t\in S$ and $\xi\in\Aa_t$ we have
			$$\tilde{E}(\xi)=\tilde{EL}(\xi)=\xi\cdot 1_{I_{1,t}},$$
			where $1_{I_{1,t}}$ is the unit in $\Mloc(I_{1,t})$, which is equal to $\Ii_{1,t}$ by Lemma~\ref{lem:intersectIdealLocalised}.
		\end{proof}
	\end{corollary}
	
	The canonical weak expectation mentioned in Corollary~\ref{cor:condExpGen} is that defined in \cite[Lemma~4.5]{BEM1}. 
	The definition of this expectation is not given in here, since under the conditions we impose this expectation is equal to the canonical local expectation.
	
	We shall now use \cite[Theorem~7.2]{BE2} to gain a saturated Fell bundle $\tilde\Ee$ over an inverse semigroup $\tilde{S}$, such that any constructed $C^*$-algebras from $\Aa$ and $\tilde\Ee$ are isomorphic via an isomorphism that entwines conditional expectations. 
	If the original action is aperiodic, the induced local action we construct will be closed and purely outer.
	
	\begin{proposition}\label{prop:saturateFellBundle}
		Let $\Ee$ be an action of $S$ on $A$. 
		There exists an inverse semigroup $\tilde{S}$ and a closed inverse semigroup action $\tilde\Ee$ of $\tilde{S}$ on $\Mloc(A)$ such that $\Mloc(A)\rtimes\tilde{S}\cong C^*(\Aa)$ via an isomorphism that maps $\Mloc(A)$ identically to itself and entwines the canonical local conditional expectations.
		Moreover, if $\Ee$ is an aperiodic action, then the action $\tilde\Ee$ is purely outer.
		
		\begin{proof}
			The saturated Fell bundle exists by \cite[Theorem~7.2]{BE2} and the bimodules associated to the action are of the form $\tilde\Ee_x=r(\Aa_{t_1})\dots r(\Aa_{t_n})\Aa_t$ for some $t_1,\dots,t_n,t\in S$, and the isomorphism described in \cite[Theorem~7.2]{BE2} maps each $\Aa_t$ identically to itself. 
			We see that the isomorphism $C^*(\tilde\Ee)\to\Mloc(A)\rtimes\tilde{S}$ must preserve the conditional expectations, since on each fibre $\Aa_x$ we have that the conditional expectation is the restriction of $\tilde{E}$ on $\Aa_t$ to a subbimodule. 
			The conditional expectation is genuine on $C^*(\Aa)$ by Corollary~\ref{cor:condExpGen}, so the action $\tilde\Ee$ is closed.
			If $\Ee$ is an aperiodic action, then each of these bimodules acts purely outerly since each $\Aa_t\cdot \Ii_{1,t}^\perp=(\Ee_t\cdot I_{1,t}^\perp)_\loc$ is purely outer by Lemma~\ref{prop:locPurOutIsAper}. 
		\end{proof}
	\end{proposition}
	
	\begin{remark}
		The inverse semigroup $\tilde{S}$ in Proposition~\ref{prop:saturateFellBundle} is called the \emph{prefix expansion} of $S$ and its construction can be found in \cite{BE2}.
		The prefix expansion of $S$ is the inverse semigroup generated by $S$ under alternative relations. 
		These relations give rise to more idempotents in $\tilde{S}$, as well as a canonical injective \emph{partial homomorphism} $\pi:S\to\tilde{S}$, which is an injective map satisfying $\pi(t)^*=\pi(t^*)$, $\pi(tu)\leq\pi(t)\pi(u)$ for all $t,u\in S$, and if $t\leq u$ then $\pi(t)\leq\pi(u)$. 
		In the context of this induced local Fell bundle, we have shown that $\tilde\Ee_t^*=\tilde\Ee_{t^*}$,  $\tilde\Ee_t\otimes_{\Mloc(A)}\tilde\Ee_u\subseteq\tilde\Ee_{tu}$ for all $t,u\in S$, and $\tilde\Ee_t\subseteq\tilde\Ee_u$ in $C^*(\tilde\Ee)$ whenever $t\leq u$, which exactly describes a partial homomorphism.
		
		The prefix expansion inverse semigroup (denoted $\mathrm{Pr}(S)$) in general differs from the refinement $\Bis(\hat{A}\rtimes S)$ in Proposition~\ref{prop-refinedAction} for two main reasons.
		The construction of the prefix inverse semigroup depends only on the inverse semigroup $S$ and is therefore independent of the $C^*$-algebra $A$ and the action $S\curvearrowright A$, whereas the refinement $\Bis(\hat{A}\rtimes S)$ can differ based on both these choices.
		For example, the idempotent lattice of $\Bis(\hat{A}\rtimes S)$ is canonically isomorphic to the ideal lattice of $A$.
		The second place where the construction of $\mathrm{Pr}(S)$ differs is a consequence of its universal property.
		Since the prefix inverse semigroup is universal for submultiplicative maps from $S$, the canonical map $S\to\mathrm{Pr}(S)$ may not be a semigroup homomorphism, but rather a submultiplicative map.
		The map $S\to \Bis(\hat{A}\rtimes S)$ however is always a semigroup homomorphism.
	\end{remark}
	
	\begin{definition}
		Let $\Ee$ be an action of $S$ on $A$. We call the action $\tilde\Ee$ of $\tilde{S}$ on $\Mloc(A)$ the \emph{induced local action}.
	\end{definition}
	
	The induced local action coming from an aperiodic action $\Ee$ gives rise to a Cartan inclusion $\Mloc(A)\subseteq\Mloc(A)\rtimes_\red S$ in the sense of Exel \cite{E1} by \cite[Theorem~4.3]{KM1}. 
	The inclusion we are interested in however is $A\subseteq A\rtimes_\ess S$, so we require a way to descend the structure of the induced inclusion to the original inclusion. 
	
	\begin{theorem}\label{thm:crossedProdEmbed}
		Let $\Ee$ be an action of $S$ on $A$.
		There is an injective homomorphism $A\rtimes_{\ess}S\hookrightarrow\Mloc(A)\rtimes_\red\tilde{S}$ that commutes with the inclusions $\Ee_t\hookrightarrow\tilde\Ee_t\subseteq\Mloc(A)\rtimes_\red\tilde{S}$ for each $t\in S$. That is, for each $t\in S$ the following diagram commutes:
		\[
			\begin{tikzpicture}[baseline=(current bounding box.west)]
			   	\node (1) at (0,1.5) {$\Ee_t$};
      			\node (2) at (2,1.5) {$A\rtimes_\ess S$};
     			\node (3) at (0,0) {$\tilde\Ee_t$};
     			\node (4) at (2,0) {$\Mloc(A)\rtimes_\red\tilde{S}$};
     			\draw[cdar] (1) -- (2);
     			\draw[cdar] (1) -- (3);
      			\draw[cdar] (2) -- (4);
     			\draw[cdar] (3) -- (4);
    		\end{tikzpicture}
		\]
		\begin{proof}
			By Proposition~\ref{prop-refinedAction} we may assume the action $\Ee$ is fine.
			Recall that $A\rtimes S$ is defined as the maximal $C^*$-completion of $A\rtimes_\alg S$, which is in turn the quotient of $\bigoplus_{t\in S}\Ee_t$ by $\Nn:=\sspan\{\xi_t\delta_t-\vartheta_{s,t}(\xi_t)\delta_s:t,s\in S, \xi_t\in\Ee_t\cdot I_{t,s}\}$. Similarly, $\Mloc(A)\rtimes\tilde{S}$ is the quotient of $\bigoplus_{t\in\tilde{S}}\tilde\Ee_t$ by the ideal $\tilde\Nn=\sspan\{\xi_t\delta_t-\tilde\vartheta_{s,t}(\xi_t)\delta_s:t,s\in \tilde{S}, \xi_t\in\tilde\Ee_t\cdot \Mloc(I_{t,s})\}$. The inclusion $\bigoplus_{t\in S}\Ee_t\to\bigoplus_{t\in\tilde{S}}\tilde\Ee_{t}$ descends to both quotients since under this inclusion as we have $\Nn\subseteq\tilde\Nn$ whereby we gain a map $A\rtimes_\alg S\to\Mloc(A)\rtimes_\alg\tilde{S}$. This map then extends to a homomorphism $i:A\rtimes S\to\Mloc(A)\rtimes S$.
			
			To see this map descends further to the essential and reduced crossed products, we first show that the map entwines conditional expectations. 
			For $\xi_t\in\Ee_t$, the element $EL(\xi_t)\in\Mloc(A)$ is defined as the multiplier in $M(I_{1,t})$ mapping $ a\mapsto \xi_t a$, where $\Ee_t\cdot I_{1,t}$ and $I_{1,t}$ are identified via $\vartheta_{1,t}$.
			By the proof of Corollary~\ref{cor:condExpGen} the element $\tilde{E}(\xi_t)\in\Mloc(A)$ is given as $\tilde{E}(\xi_t)=\xi_t\cdot 1_{I_{1,t}}$.
			Thus for $a\in I_{1,t}$ we have $\tilde{E}(i(\xi_t))a=\xi_t\cdot a$, where $\tilde\Ee_t\cdot\Ii_{1,t}$ is identified with $\Ii_{1,t}$ by applying Proposition~\ref{prop:tensorProdLocalMap} to the isomorphism $\vartheta_{1,t}$.
			Since the bimodules $\Ee_t$ span a dense subspace of $A\rtimes_\ess S$, we see that $EL=\tilde{E}\circ i$. 
			Since $A\rtimes_\ess S$ is the quotient of $A\rtimes S$ by the ideal $N_{EL}=\{a\in A\rtimes S:EL(a^*a)=0\}=\{a\in A\rtimes S:\tilde{E}\circ i(a^*a)=0\}$ we have $i(N_{EL})\subseteq N_{\tilde{E}}$, so $i$ descends to a map on $A\rtimes_\ess S$. 
			
			Lastly we show that $i$ is injective on $A\rtimes_\ess S$. To see this, we note that an element $a\in A\rtimes_\ess S$ is mapped to zero under $i$ if and only if it satisfies $\tilde{E}(i(a^*a))=EL(a^*a)=0$, giving $a=0$.
		\end{proof}
	\end{theorem}
	
	Theorem~\ref{thm:crossedProdEmbed} allows the dynamic and algebraic structures of the actions $\Ee$ and $\tilde\Ee$ to be encoded in the one algebra $\Mloc(A)\rtimes_\red\tilde{S}$, and provides a setting useful for computations. 
	Throughout the rest of this article we identify the modules $\Ee_t,\tilde\Ee_t$, and the algebras $A,\Mloc(A)$ and $A\rtimes_\ess S$ with their images in $\Mloc(A)\rtimes_\red\tilde{S}$ via Theorem~\ref{thm:crossedProdEmbed}.
	
	\begin{remark}
		The refinement from $S$ to $\Bis(\hat{A}\rtimes S)$ is required since taking local multiplier algebras does not commute with taking inductive limits.
		Given a countable collection $(A_n)_{n\in\NN}$ of non-zero $C^*$-algebras, one may consider the inductive system given by algebras $B_n=\bigoplus_{j=1}^n A_j$ with maps $\iota_n=(\Id_{B_n},0):B_n\to B_n\oplus A_{n+1}=B_{n+1}$.
		The inductive limit of this system is the direct sum $\bigoplus_{n\in\NN} A_n$.
		Finite direct sums agree with finite direct products, so Lemma~\ref{lem:locMultAddProd} gives that $\Mloc(\cdot)$ preserves finite direct sums.
		By Lemma~\ref{lem:locMultAddProd}, the local multipiler algebra of the inductive limit of the system $B_n$ is then 
		$$\Mloc\left(\varinjlim B_n\right)=\Mloc\left(\bigoplus_{n\in\NN}A_n\right)=\prod_{n\in\NN}\Mloc(A_n),$$
		which is a unital $C^*$-algebra.
		The inductive limit of local multiplier algebras is however
		$$\varinjlim \Mloc(B_n)=\varinjlim \bigoplus_{j=1}^n\Mloc(A_j)=\bigoplus_{n\in\NN}\Mloc(A_n).$$
		This cannot be isomorphic to the local multiplier algebra of the inductive limit, since direct sums of infinitely many non-zero $C^*$-algebras are never unital.
	\end{remark}

	\section{Aperiodic dynamical inclusions}
	
	With Theorem~\ref{thm:crossedProdEmbed} we see that $A\rtimes_\ess S$ embeds into $\Mloc(A)\rtimes_\red\tilde{S}$ in a way that preserves the inclusion $A\subseteq\Mloc(A)$.
	As briefly mentioned, the inclusion $A\subseteq\Mloc(A)$ is problematic, and so we wish to ensure that the intersection of $A\rtimes_\ess S$ and $\Mloc(A)$ in $\Mloc(A)\rtimes_\red\tilde{S}$ is as small as possible. 
	The local conditional expectation $EL:A\rtimes_\ess S\to\Mloc(A)$ restricts to a ${}^*$-homomorphism on at least $A$, since $EL|_A=\Id_A$, and so to minimise the intersection between $\Mloc(A)$ and $A\rtimes_\ess S$ we investigate the circumstances under which $EL$ restricts to a ${}^*$-homomorphism only on $A$.
	
	\subsection{Inclusions with minimal multiplicative domain}
	
	\begin{definition}
		Let $A\subseteq B$ be an inclusion of $C^*$-algebras and let $E:B\to\tilde{A}$ be a generalised expectation. 
		We say that $E$ has \emph{minimal multiplicative domain} (MMD) if the following subset is equal to $A$:
		$$\mu(E):=\{b\in B: E(b^*b)=E(b^*)E(b), E(bb^*)=E(b)E(b^*)\}.$$
	\end{definition}
	
	The set $\mu(E)$ for an expectation $E$ is the (two-sided) multiplicative domain of $E$; the largest $C^*$-subalgebra of $B$, on which $E$ restricts to a ${}^*$-homomorphism. This characterisation is adapted from Choi \cite[Theorem~3.1]{C1} in which the multiplicative domain for 2-positive maps between $C^*$-algebras is defined. This formulation applies here as generalised expectations are completely positive by definition, and the two-sided condition ensures that $\mu(E)$ is ${}^*$-closed. The `minimal' descriptor of the definition of (MMD) stems from the fact that we always have $A\subseteq\mu(E)$, as $E$ restricts to the identity on $A$. 
	
	If $E$ is faithful, then $E$ restricts to an injective ${}^*$-homomorphism $\mu(E)\hookrightarrow\tilde{A}$. Thus, in the case of a faithful genuine conditional expectation $E:B\to A$, the multiplicative domain of $E$ is always minimal.
	
	\begin{lemma}
		Let $\Ee$ be an action of $S$ on $A$ and let $\tilde\Ee$ be its induced local action. Consider $A\rtimes_\ess S$ as a subalgebra of $\Mloc(A)\rtimes_\red\tilde{S}$ as in Theorem~\textup{\ref{thm:crossedProdEmbed}}. Then the multiplicative domain of $EL:A\rtimes_\ess S\to\Mloc(A)$ is equal to the intersection $\Mloc(A)\cap A\rtimes_\ess S$, where the intersection is taken in $\Mloc(A)\rtimes_\red\tilde{S}$.
		\begin{proof}
			Let $\tilde{E}:\Mloc(A)\rtimes_\red\tilde{S}\to\Mloc(A)$ be the canonical conditional expectation for the action $\tilde\Ee$. 
			By Theorem~\ref{thm:crossedProdEmbed} the inclusion $A\rtimes_\ess S\subseteq\Mloc(A)\rtimes_\red\tilde{S}$ entwines conditional expectations, so we have $\tilde{E}|_{A\rtimes_\ess S}=EL$.
			Since $\tilde{E}$ is a faithful genuine expectation we have $\mu(\tilde{E})=\Mloc(A)$. Thus
			\begin{align*}
			\mu(EL)=\mu(\tilde{E}|_{A\rtimes_\ess S})&=\mu(\tilde{E})\cap A\rtimes_\ess S=\Mloc(A)\cap A\rtimes_\ess S.\qedhere
			\end{align*}
		\end{proof}
	\end{lemma}
	
	\begin{corollary}\label{cor:mmdIntersectProp}
		The local expectation $EL$ for the inclusion $A\subseteq A\rtimes_\ess S$ has minimal multiplicative domain if and only if $\Mloc(A)\cap A\rtimes_\ess S=A$.
	\end{corollary}
	
	Corollary~\ref{cor:mmdIntersectProp} does not work generally for the modules $\Ee_t\subseteq A\rtimes_\ess S$: it fails even for non-unital ideals $I\triangleleft A$ since Lemma~\ref{lem:AintersectMlocIdeal} gives $\Mloc(I)\cap A\rtimes_\ess S= M(I)\cap A\supseteq I$. However, we always will have containment $\Ee_t\subseteq\tilde\Ee_t\cap A\rtimes_\ess S$.
	
	\begin{lemma}\label{lem:recoverActionModuleLocal}
		Let $\Ee$ be an action such that $EL:A\rtimes_\ess S\to\Mloc(A)$ has minimal multiplicative domain. For $t\in S$ we have $\Ee_t=(\tilde\Ee_t\cap A\rtimes_\ess S)\cdot s(\Ee_t)$.
		\begin{proof}
			The inclusion $\Ee_t\subseteq (\tilde\Ee_t\cap A\rtimes_\ess S)\cdot s(\Ee_t)$ follows since $\Ee_t$ is contained in both $\tilde\Ee_t$ and $A\rtimes_\ess S$, and $\Ee_t=\Ee_t\cdot s(\Ee_t)$. 
			We shall show that $\tilde\Ee_t\cap A\rtimes_\ess S$ is a Hilbert $A$-bimodule, and then the reverse inclusion follows as the source ideals will be equal.
						
			For $\xi,\eta\in\tilde\Ee_t\cap A\rtimes_\ess S$, both $\langle\xi,\eta\rangle=\xi^*\eta$ and $\llangle\xi,\eta\rrangle=\xi\eta^*$ belong to $\Mloc(A)$ as $\tilde\Ee_t$ is a Hilbert $\Mloc(A)$-bimodule, and belong to $A\rtimes_\ess S$ as $\xi,\eta\in A\rtimes_\ess S$ which is closed under its own multiplication. Thus the inner products of $\tilde\Ee_t\cap A\rtimes_\ess S$ take values in $\Mloc(A)\cap A\rtimes_\ess S$, which is $A$ by Corollary~\ref{cor:mmdIntersectProp}. Then $\tilde\Ee_t\cap A\rtimes_\ess S$ is closed under the left and right $A$-multiplications as both $\tilde\Ee_t$ and $A\rtimes_\ess S$ are, and is norm-closed as an intersection of closed subsets of $\Mloc(A)\rtimes_\red\tilde{S}$.
		\end{proof}
	\end{lemma}
	
	\subsection{Slice reconstruction}
	
	With these stronger results for inclusions with minimal multiplicative domain we can now define the class of inclusions we wish to study.
	
	\begin{definition}
		Let $A\subseteq B$ be an inclusion of $C^*$-algebras.
		We say the inclusion is an \emph{aperiodic dynamical inclusion} if the following conditions hold:
		\begin{enumerate}
			\item $A\subseteq B$ is non-degenerate;
			\item there exists an inverse subsemigroup $S\subseteq\Ss(A,B)$ such that the tautological action of $S$ on $A$ is aperiodic; and
			\item there exists a faithful pseudo-expectation $E:B\to I(A)$ with minimal multiplicative domain.
		\end{enumerate}
	\end{definition}

	\begin{proposition}\label{prop-aperDynIncClassification}
		Let $A\subseteq B$ be an aperiodic dynamical inclusion.
		Let $E:B\to I(A)$ be the associated faithful pseudo-expectation with minimal multiplicative domain and let $S\subseteq \Ss(A,B)$ be the distinguished inverse subsemigroup acting aperiodically.
		Then there is an isomorphism $\varphi:A\rtimes_\ess S\to B$ that entwines conditional expectations and hence restricts to the identity map on $A$.
		\begin{proof}
			Aperiodic dynamical inclusions satisfy the hypotheses of Theorem~\ref{thm:mainAperIso}, which then gives the desired results.
		\end{proof}
	\end{proposition}
	
	\begin{corollary}
		Let $\Ee$ be an aperiodic action such that the local expectation $EL:A\rtimes_\ess S\to \Mloc(A)$ has minimal multiplicative domain.
		Then $A\subseteq A\rtimes_\ess S$ is an aperiodic dynamical inclusion, and up to isomorphism every aperiodic dynamical inclusion is of this form.
	\end{corollary}
	
	One of the statements in \cite[Theorem~5.6]{KM1} is that the bisection inverse semigroup $\Bis(\hat{A}\rtimes S)$ for a closed and purely outer action is isomorphic to the slice inverse semigroup $\Ss(A,A\rtimes_\red S)$. Slices for the inclusion $A\subseteq A\rtimes_\red S$  can then be recovered from their intersections with the bimodules $\Ee_t$.
	Thus all non-trivial slices for the inclusion $A\subseteq A\rtimes_\red S$ inherit pure-outerness, and there can be no wayward slices outside of those coming from the action.
	We shall apply this to the inclusion $\Mloc(A)\subseteq\Mloc(A)\rtimes_\red\tilde{S}$, and then use (MMD) to descend properties to the inclusion $A\subseteq A\rtimes_\ess S$.
	
	\begin{lemma}\label{lem:sliceReconstruction}
		Let $\Ee:S\curvearrowright A$ be a closed and purely outer action. Let $X\subseteq A\rtimes_\red S$ be a slice. Then $X=\overline{\sum_{t\in S}X\cap\Ee_t}$.
		\begin{proof}
			The inclusion $\overline{\sum_{t\in S}X\cap\Ee_t}\subseteq X$ is clear. The proof of \cite[Theorem~5.6]{KM1} shows that $X$ is the closure of the span of slices $\Ee_t\cdot J_t$ for some $t\in T\subseteq S$ and ideals $J_t\triangleleft A$. Each $\Ee_t\cdot J_t$ is contained in $\Ee_t$ and $X$, so we gain the desired result.
		\end{proof}
	\end{lemma}
	
	In the remainder of this section we shall fix an action $\Ee$ of $S$ on $A$ such that the canonical local expectation $EL:A\rtimes_\ess S\to\Mloc(A)$ has minimal multiplicative domain.
	We shall denote by $\tilde\Ee$ the induced local action of $\tilde{S}$ on $\Mloc(A)$ arising from $\Ee$.
	By Proposition~\ref{prop-refinedAction} we may and will without loss of generality assume the action $\Ee$ to be fine.
	
	\begin{corollary}\label{cor:essCrossDetectsSlices}
		Let $\tilde{Y}\subseteq\Mloc(A)\rtimes_\red\tilde{S}$ be a slice for the induced inclusion $\Mloc(A)\subseteq \Mloc(A)\rtimes_\red\tilde{S}$. Then $\tilde{Y}=\{0\}$ if and only if $\tilde{Y}\cap A\rtimes_\ess S=\{0\}$, where the intersection is taken in $\Mloc(A)\rtimes_\red\tilde{S}$ using Theorem~\textup{\ref{thm:crossedProdEmbed}}.
		\begin{proof}
			If $\tilde{Y}=\{0\}$ then clearly $\tilde{Y}\cap A\rtimes_\ess S=\{0\}$. 
			Conversely if $\tilde{Y}\neq\{0\}$ then Lemma~\ref{lem:sliceReconstruction} implies that for some $t\in S$ we have $\tilde{Y}\cap\tilde\Ee_t\neq\{0\}$, which is a subbimodule of $\tilde\Ee_t$. Lemma~\ref{lem:detectSubbimodsOfLoc} then gives $\{0\}\neq \tilde{Y}\cap\tilde\Ee_t\cap\Ee_t\subseteq\tilde{Y}\cap A\rtimes_\ess S$. 
		\end{proof}
	\end{corollary}
	
	If $X\subseteq A\rtimes_\ess S$ is a slice, one would want to analyse a corresponding slice $\breve{X}\subseteq\Mloc(A)\rtimes_\red\tilde{S}$ generated by $X$.
	
	\begin{lemma}
		Let $X\subseteq A\rtimes_\ess S$ be a slice for the inclusion $A\subseteq A\rtimes_\ess S$. Considering $A\rtimes_\ess S\subseteq\Mloc(A)\rtimes_\red\tilde{S}$ as in Theorem~\textup{\ref{thm:crossedProdEmbed}}, $\breve{X}:=\overline\sspan\Mloc(A)\cdot X\cdot\Mloc(A)$ is a slice for $\Mloc(A)\subseteq\Mloc(A)\rtimes_\red\tilde{S}$ and satisfies $X=(\breve{X}\cap A\rtimes_\ess S)\cdot s(X)$.
		\begin{proof}
			We first show that $\breve{X}$ is a slice for $\Mloc(A)\subseteq\Mloc(A)\rtimes_\red\tilde{S}$. Fix $x,y\in X$ and $a\in M(J)$ for some $J\in\Ii_e(A)$. Let $I:= s(J\cdot X)\oplus s(J\cdot X)^\perp$. This ideal is essential in $A$ by construction, and satisfies $X\cdot I=J\cdot X$.
			By the Cohen-Hewitt Factorisation Theorem for all $b\in I$ there exists $z\in X$ and $c\in J$ such that $yb=cz$. We then have $(x^*ay)b=x^*(ac)z\in X^*M(J)JX=X^*JX=s(J\cdot X)\subseteq I$. 
			We also note that $b(x^*ay)=(xb^*)^*ay$ belongs to  $(X\cdot I)^*M(J)X=(J\cdot X)^*M(J)X=X^*JM(J)X\subseteq I$.
			Thus $x^*ay$ is a multiplier on $I$, so belongs to the local multiplier algebra. Since $||x^*ay||\leq ||x||\cdot||y||\cdot||a||$, it follows that $x^*\Mloc(A)y\subseteq\Mloc(A)$ for all $x,y\in X$. Thus $\breve{X}^*\Mloc(A)\breve{X}=\cspan\Mloc(A)\cdot X^*\cdot\Mloc(A)\cdot X\cdot\Mloc(A)\subseteq\Mloc(A)^3=\Mloc(A)$.
			Symmetrically $\breve{X}\Mloc(A)\breve{X}^*\subseteq\Mloc(A)$, so $\breve{X}$ is a slice.
			
			Similarly to Lemma~\ref{lem:recoverActionModuleLocal}, the inclusion $X\subseteq (\breve{X}\cap A\rtimes_\ess S)\cdot s(X)$ is clear. For the reverse inclusion, we see that $\breve{X}\cap A\rtimes_\ess S$ has inner products taking value in $\Mloc(A)\cap A\rtimes_\ess S$, which is equal to $A$ by Corollary~\ref{cor:mmdIntersectProp}. By cutting down with the source ideal $s(X)$, we then gain the desired equality.
		\end{proof}
	\end{lemma}

	\begin{corollary}\label{cor:sliceDetectsInGend}
	Let $X\subseteq A\rtimes_\ess S$ be a slice and let $\breve{X}=\cspan\Mloc(A)\cdot X\cdot\Mloc(A)$ be the slice in $\Mloc(A)\rtimes_\red\tilde{S}$ generated by $X$. 
	Then for any Hilbert $A$-subbimodule $Y\subseteq \breve{X}\cap A\rtimes_\ess S$ we have $Y=\{0\}$ if and only if $Y\cap X=\{0\}$.
		\begin{proof}
			The source ideal of $\breve{X}$ annihilates $\Mloc(s(X)^\perp)$ in $\Mloc(A)$ since 
			$$\breve{X} s(X)^\perp\subseteq \cspan\Mloc(A)X\Mloc(s(X)^\perp)=\{0\}.$$
			Thus $s(\breve{X})\subseteq\Mloc(s(X)^\perp)^\perp=\Mloc(s(X))$. We then have $s(\breve{X}\cap A\rtimes_\ess S)\subseteq\Mloc(s(X))\cap A\rtimes_\ess S$, which is equal to $M(s(X))\cap A$ by Lemma~\ref{lem:AintersectMlocIdeal} and Corollary~\ref{cor:mmdIntersectProp}. 

			We claim that $s(X)$ is an essential ideal of $s(\breve{X}\cap A\rtimes_\ess S)$.
			Note that the only element of $\Mloc(s(X))$ that annihilates $s(X)$ is $0$, and so the annihilator of $s(X)$ in $\Mloc(A)$ is $\Mloc(s(X)^\perp)$ by Lemma~\ref{lem:locMultAddProd}.
			We then see that $X\cdot\Mloc(A)\cdot s(X)^\perp\subseteq X\cdot\Mloc(s(X)^\perp)=\{0\}$, giving $\breve{X}\cdot s(X)^\perp=\{0\}$.
			Thus $s(X)$ must be an essential ideal of $s(\breve{X}\cap A\rtimes_\ess S)$.
			Hence if $Y\subseteq\breve{X}\cap A\rtimes_\ess S$ is a non-zero Hilbert $A$-subbimodule, then $Y\cdot s(X)\neq \{0\}$ as $s(X)$ is essential in $s(\breve{X}\cap A\rtimes_\ess S)$. 
			
			Thus $Y\cdot s(X)\subseteq(\breve{X}\cap A\rtimes_\ess S)\cdot s(X)=X$, so $\{0\}\neq Y\cdot s(X)\subseteq X\cap Y$. Contrapositively, $Y\cap X=\{0\}$ gives $Y=\{0\}$.
		\end{proof}
	\end{corollary}	
	
	The conditional expectation having minimal multiplicative domain allows us to analyse slices of $A\subseteq A\rtimes_\ess S$ by analysing slices of $\Mloc(A)\subseteq\Mloc(A)\rtimes_\red\tilde{S}$ and descending down to the smaller inclusion.
	If the action $\Ee$ is aperiodic, then the induced action $\tilde\Ee$ is closed and purely outer by Proposition~\ref{prop:saturateFellBundle}, and so the inclusion $\Mloc(A)\subseteq\Mloc(A)\rtimes_\red\tilde{S}$ satisfies the hypotheses of Lemma~\ref{lem:sliceReconstruction}. 
	We use this now to show that all slices in an aperiodic dynamical inclusion act aperiodically, and all choices of inverse subsemigroup $S\subseteq\Ss(A,B)$ that densely span $B$ give rise to the same essential crossed product.
	From this point on we assume that the action $\Ee$ is aperiodic, so that the induced local action $\tilde\Ee$ is purely outer.
	
	\begin{lemma}\label{lem:aperActionDetectsSlices}
		Suppose that the action $\Ee$ is aperiodic.
		Let $X\subseteq A\rtimes_\ess S$ be a slice for the inclusion $A\subseteq A\rtimes_\ess S$. Then $X=\{0\}$ if and only if $X\cap\Ee_t=\{0\}$ for all $t\in S$.
		\begin{proof}
			The `only if' direction is clear. 
			
			Suppose that $X\cap\Ee_t=\{0\}$ for all $t\in S$. Let $\breve{X}:=\overline\sspan\Mloc(A)\cdot X\cdot\Mloc(A)$ be the slice for the larger inclusion $\Mloc(A)\subseteq\Mloc(A)\rtimes_\red\tilde{S}$ generated by $X$. The equality $X\cap\Ee_t=\{0\}$ implies that $X\cap(\Ee_t\cap\breve{X}\cap A\rtimes_\ess S)=\{0\}$ for each $t\in S$, and so by Corollary~\ref{cor:sliceDetectsInGend} we have $\{0\}=\Ee_t\cap\breve{X}\cap A\rtimes_\ess S=\Ee_t\cap\breve{X}\supseteq \Ee_t\cap\breve{X}\cap\tilde\Ee_t$. The module $\tilde\Ee_t\cap\breve{X}$ is a subbimodule of $\tilde\Ee_t$ and $\Ee_t$ detects subbimodules of its local multiplier module by Lemma~\ref{lem:detectSubbimodsOfLoc}.
			Thus $\Ee_t\cap\breve{X}=\{0\}$ implies that $\tilde\Ee_t\cap\breve{X}=\{0\}$.
			Since the action $\Ee$ is aperiodic, the induced local action is closed and purely outer by Proposition~\ref{prop:saturateFellBundle}, and so satisfies the criteria of Lemma~\ref{lem:sliceReconstruction}.
			Hence $\breve{X}\cap\tilde\Ee_t=\{0\}$ for all $t\in S$ if and only if $\breve{X}=\{0\}$.
		\end{proof}
	\end{lemma}
	
	\begin{corollary}\label{cor:aperActionEssentialSubslice}
		Suppose that the action $\Ee$ is aperiodic.
		For $X\in\Ss(A,A\rtimes_\ess S)$, the subslice $X_\ess:=\overline{\sum_{t\in S}X\cap\Ee_t}$ has zero orthogonal complement in $X$, and so has non-zero intersection with all non-zero subslices of $X$.
		Equivalently, $s(X_\ess)$ is an essential ideal of $s(X)$.
		\begin{proof}
			We see $X_\ess^\perp\cap X$ satisfies $X_\ess^\perp\cap X\cap\Ee_t=\bigcap_{u\in S}(X\cap\Ee_u)^\perp\cap(X\cap\Ee_t)\subseteq(X\cap\Ee_t)^\perp\cap (X\cap\Ee_t)=\{0\}$.
			Hence $X_\ess^\perp\cap X=\{0\}$ by Lemma~\ref{lem:aperActionDetectsSlices}. If $Y\subseteq X$ is a non-zero subslice, then $Y\cap X_\ess=Y\cap\overline{\sum_{t\in S}X\cap\Ee_t}\supseteq \overline{\sum_{t\in S}Y\cap\Ee_t}$, which is zero if and only if $Y=\{0\}$ by Lemma~\ref{lem:aperActionDetectsSlices}.
			That this is equivalent to $s(X_\ess)\triangleleft s(X)$ being an essential ideal follows from the fact that any subbimodule $Y\subseteq X$ satisfies $Y=X\cdot s(Y)$.
		\end{proof}
	\end{corollary}
	
	\begin{theorem}\label{thm-mainAperDynIncProperty}
		Let $A\subseteq B$ be an aperiodic dynamical inclusion.
		Then the tautological action of $\Ss(A,B)$ on $A$ is aperiodic.
		If $T\subseteq\Ss(A,B)$ is an inverse subsemigroup that spans a dense subalgebra of $B$, then there are isomorphisms
		$$B\cong A\rtimes_\ess T\cong A\rtimes_\ess\Ss(A,B).$$
		Moreover these isomorphisms map $A$ identically to itself and entwine conditional expectations.
		\begin{proof}
			Let $S\subseteq\Ss(A,B)$ be an inverse subsemigroup that acts aperiodically and spans a dense subalgebra of $B$ (note that such a semigroup exists since $A\subseteq B$ is an aperiodic dynamical inclusion).
			By Theorem~\ref{thm:mainAperIso} there is an isomorphism $B\cong A\rtimes_\ess S$ that maps $A$ identically to itself and entwines conditional expectations.
			
			Let $X\subseteq A\rtimes_\ess S\cong B$ be a slice and let $Y:=X\cdot (X\cap A)^\perp$.
			By Corollary~\ref{cor:aperActionEssentialSubslice} the subslice $Y_\ess:=\overline{\sum_{t\in S}Y\cap \Ee_t}$ has zero orthogonal complement in $Y$ and so $s(Y_\ess)$ is an essential ideal of $s(Y)$.
			For each $t\in S$ the ideal $s(Y\cap\Ee_t\cdot I_{1,t}^\perp)$ is an essential ideal of $s(Y\cap\Ee_t)$ since $Y\cap A=\{0\}$. 
			Applying \cite[Lemma~5.12]{KM2} shows that $Y\cap\Ee_t$ is an aperiodic $A$-bimodule since $Y\cap\Ee_t\cdot I_{1,t}^\perp$ is.
			Thus $Y$ is aperiodic as the closed linear span of aperiodic subbimodules $Y\cap\Ee_t$, and so $X$ acts aperiodically on $A$.
			Thus any choice of $T\subseteq\Ss(A,B)$ makes the inclusion $A\subseteq B$ an aperiodic dynamical inclusion, and so $B\cong A\rtimes_\ess T$ by Theorem~\ref{thm:mainAperIso}.
			Particularly the choice $T=\Ss(A,B)$ gives $B\cong A\rtimes_\ess\Ss(A,B)$.
		\end{proof}
	\end{theorem}
	
	The topology on the dual groupoid $\hat{A}\rtimes\Ss(A,B)$ has a basis given by slices of the inclusion $A\subseteq B$.
	Lemma~\ref{lem:aperActionDetectsSlices} shows that any non-zero slice $X\in\Ss(A,B)$ intersects at least one $\Ee_t$ for some $t\in S$.
	In the topology of the groupoid, we see then that any open subset of $\hat{A}\rtimes\Ss(A,B)$ must intersect the open bisection defined by $\Ee_t$. 
	Thus the bisections defined by $\Ee_t$ for $t\in S$ cover a dense subset of the dual groupoid for the full slice action.
	
	\begin{corollary}\label{cor:grpdDenseRange}
		Suppose the action $\Ee$ is aperiodic.
		Let $B=A\rtimes_\ess S$. 
		The canonical groupoid homomorphism $\phi:\hat{A}\rtimes S\to\hat{A}\rtimes\Ss(A,B)$, $\phi[t,[\pi]]=[\Ee_t,[\pi]]$ has open range.
		If for each $t\in S$ and each ideal $I\triangleleft A$ there is an idempotent $e\in S$ with $\Ee_t\cdot I=\Ee_e$, then $\phi$ is injective.
		\begin{proof}
			Each $\phi[t,[\pi]]$ belongs to $\hat\Ee_t$ for each $[t,[\pi]]\in A\rtimes S$ and $\bigcup_{t\in S}\hat\Ee_t$ is open in $A\rtimes\Ss(A,B)$ as each $\hat\Ee_t$ is an open bisection, and so the image of $\phi$ is open.
			
			The last part giving that $\phi$ is injective follows from Lemma~\ref{lem-enoughIdemsGivesInjGrpdMap}.
		\end{proof}
	\end{corollary}

	\begin{remark}
		The map $\phi$ in Corollary~\ref{cor:grpdDenseRange} is not necessarily injective, since the semigroup $\Ss(A,B)$ may have a more refined lattice of idempotents than $S$.
		For example if $S$ is a group (viewed as an inverse semigroup) acting on a $C^*$-algebra $A$, then the germ relation is trivial since the only idempotent in a group is the identity.
		The inverse semigroup $\Ss(A,B)$ will always contain the ideal lattice of $A$ as idempotents, and so the quotient may identify more germs.
		The condition of Lemma~\ref{lem-enoughIdemsGivesInjGrpdMap} ensures there are enough idempotents, which holds if in particular the action $\Ee$ is fine.
		
		We do not know if the map $\phi$ is automatically injective for aperiodic inverse semigroup actions.
		If the action is topologically non-trivial, that is, the bimodules $\Ee_t\cdot I_{1,t}$ induce partial homeorphisms that do not fix any open subset of $\hat{A}$, then $\phi$ is injective.
		Every topologically non-trivial action is aperiodic by \cite[Theorem~8.1]{KM3}, but the converse is not known (unless the algebra $A$ contains an essential ideal that is simple or of Type I).
		If there exists a $C^*$-algebra $A$ and a non-zero Hilbert $A$-bimodule $X$ that is aperiodic but not topologically non-trivial, then $X$ would generate such an action on $A$ where the map $\phi$ would fail to be injective.
	\end{remark}


\begin{thebibliography}{00}
		
		\bibitem{AM1} P. Ara, M. Matthieu, \emph{Local Multipliers of $C^*$-algebras}, Springer Monographs in Mathematics, Springer London (2003), https://doi.org/10.1007/978-1-4471-0045-4
		
		\bibitem{BE1} A. Buss, E. Exel, \emph{Fell bundles over inverse semigroups and twisted {\'e}tale groupoids}, Journal of Operator Theory (2012), Vol. 67, No. 1 (Winter 2012), pp. 153-205 (53 pages)
		
		\bibitem{BE2} A. Buss, E. Exel, Alcides Buss. Ruy Exel. \emph{Inverse semigroup expansions and their actions on $C^*$-algebras}, Illinois J. Math. 56 (4) 1185 -- 1212, Winter 2012. https://doi.org/10.1215/ijm/1399395828
		
		\bibitem{BEM1} A. Buss, R. Exel, R. Meyer, \emph{Reduced {$C^*$}-algebras of Fell bundles over inverse semigroups}, Israel J. Math. 220 no.1 (2017), pp. 225 -- 274
		
		\bibitem{BM1} A. Buss, R. Meyer, \emph{Inverse Semigroup Actions on Groupoids}, Rocky Mountain J. Math. 47 (1) 53 - 159, 2017. https://doi.org/10.1216/RMJ-2017-47-1-53
		
		\bibitem{C1} M. D. Choi, \emph{A Schwarz inequality for positive linear maps on {$C^*$}-algebras}, Illinois J. Math. 18 (1974), 565--574. 
		
		
		\bibitem{E1} R. Exel, \emph{Noncommutative Cartan sub-algebras of {$C^*$}-algebras}, The New York Journal of Mathematics [electronic only] 17 (2011): 331-382. <http://eudml.org/doc/223366>.
		
		\bibitem{E2} R. Exel, \emph{Non-Hausdorff {\'e}tale groupoids}, Proceedings of the American Mathematical Society 139(3) DOI:10.1090/S0002-9939-2010-10477-X
		
		\bibitem{F1} M. Frank, \emph{Injective envelopes and local multipliers of {$C^*$}-algebras}, Int. Math. J. 1(2002), no. 6, 611 -- 620 https://doi.org/10.48550/arXiv.math/9910109
		
		\bibitem{Green1} P. Green, \emph{The local structure of twisted covariance algebras}, Acta Math. 140 191 -- 250, 1978. https://doi.org/10.1007/BF02392308
		
		\bibitem{G1} H. Gonshor, \emph{Injective hulls of {$C^*$}-algebras. II}, Proceedings of the American Mathematical Society, Vol. 24, No. 3 (Mar., 1970), pp. 486 -- 491 
		
		\bibitem{H1} M. Hamana, Injective envelopes of {$C^*$}-algebras, J. Math. Soc. Japan 31 (1979), no. 1, 181–197, doi: 10.2969/jmsj/03110181. MR 519044.
		
		\bibitem{K1} A. Kumjian, \emph{On {$C^*$}-diagonals}, Canadian Journal of Mathematics, Volume 38, Issue 4 , 01 August 1986 , pp. 969 - 1008 DOI: https://doi.org/10.4153/CJM-1986-048-0
		
		\bibitem{KM3} B. K. Kwa\'sniewski, R. Meyer, \emph{Aperiodicity, topological freeness and pure outerness: from group actions to Fell bundles}, Studia Math. 241 (2018), 257-302  https://doi.org/10.48550/arXiv.1611.06954
		
		\bibitem{KM1} B. K. Kwa\'sniewski, R. Meyer, \emph{On Noncommutative Cartan {$C^*$}-subalgebras}, Trans. Amer. Math. Soc. 373 (2020), p. 8697-8724 https://doi.org/10.48550/arXiv.1908.07217
		
		\bibitem{KM2} B. K. Kwa\'sniewski, R. Meyer, \emph{Essential Crossed Products for Inverse Semigroup Actions: Simplicity and Pure Infiniteness}, Doc. Math. 26 (2021), 271-335 https://doi.org/10.48550/arXiv.1906.06202
		
		\bibitem{KM4} B. K. Kwa\'sniewski, R. Meyer, \emph{Aperiodicity, the almost extension property and uniqueness of pseudo-expectations}, International Mathematics Research Notices (IMRN), published online 07 June 2021 https://doi.org/10.48550/arXiv.2007.05409
		
		\bibitem{Lance1} E. C. Lance, \emph{Hilbert {$C^*$}-Modules: A Toolkit for Operator Algebraists}, London Mathematical Society Lecture Note Series, vol. \textbf{210}, Cambridge University Press, Cambridge, 1995. https://doi.org/10.1017/CBO9780511526206
		
		\bibitem{L1} X. Li, \emph{Every classifiable simple {$C^*$}-algebra has a Cartan subalgebra}, Invent. math. 219, 653–699 (2020). https://doi.org/10.1007/s00222-019-00914-0
		
		\bibitem{P1} D. Pitts, \emph{Normalizers and Approximate Units for Inclusions of {$C^*$}-Algebras}, (arXiv:2109.00856 [math.OA]), (2022), (To appear in Indiana University Mathematics Journal) https://doi.org/10.48550/arXiv.2109.00856
		
		\bibitem{R2} J. Renault, \emph{A groupoid approach to {$C^*$}-algebras}, Lecture Notes in Mathematics, Vol. \textbf{793} Springer-Verlag Berlin, Heidelberg, New York, 1980. https://doi.org/10.1007/BFb0091072
		
		\bibitem{R1} J. Renault, \emph{Cartan subalgebras in {$C^*$}-algebras}, Irish Math. Soc. Bull. \textbf{61} (2008), 29--63. DOI:10.33232/BIMS.0061.29.63
		
		\bibitem{Raad1} A. Imad Raad, \emph{A Generalization of Renault's Theorem for Cartan Subalgebras}, (arXiv:2101.03265 [math.OA]), (2021), (Accepted for publication in Proc. Amer. Math. Soc) https://doi.org/10.48550/arXiv.2101.03265
		
		\bibitem{T1} J. Taylor, \emph{Aperiodic dynamical inclusions of $C^*$-algebras}, Georg-August-Universit\"at G\"ottingen (2022), http://dx.doi.org/10.53846/goediss-9727
	\end{thebibliography}
\end{document}